\documentclass[a4paper,11pt,openany]{article}
\usepackage{indentfirst}
\usepackage[greek,french,english]{babel}
\usepackage[iso-8859-7]{inputenc}

\usepackage{color,listings,titlesec}
\usepackage[headings]{fullpage}
\usepackage{amsfonts}
\usepackage{fancyhdr}
\usepackage{graphicx}
\usepackage{amsmath}
\usepackage{marvosym}
\usepackage{hyperref}
\newcommand{\g}{\mathfrak{g}}
\newcommand{\h}{\mathfrak{h}}

\newcommand{\ad}{\mathrm{ad}}

\title{An isomorphism between the reduction algebra and the invariant differential operators on Lie groups.}
\author{P. Batakidis\footnote{Department of Mathematics \& Statistics, University of Cyprus, Cyprus. E-mail: batakid@ucy.ac.cy}}
\begin{document}
\maketitle
\begin{abstract}
\noindent Let $G$ be a connected and simply connected Lie group with Lie algebra $\g$ of finite dimension. Let $\h\subset\g$ be a subalgebra, $\lambda$ a character of $\h$, $\rho(H)=-\frac{1}{2}\mathrm{Tr}\;\mathrm{ad}_\g(H)$, $H\in\h$ and $\epsilon$ a formal parameter. We prove an algebra isomorphism $H^0_{(\epsilon)}(\h_{\lambda}^{\bot},d^{(\epsilon)}_{\h_{\lambda}^{\bot},\mathfrak{q}})\simeq(U_{(\epsilon)}(\g)/U_{(\epsilon)}(\g)\h_{\lambda+\rho})^{\h}$ where the left algebra is the reduction algebra over the affine space $-\lambda+\h^{\bot}$, and at the right part $U_{(\epsilon)}(\g)$ is a deformation of the universal enveloping algebra $U(\g)$ of $\g$ and $U_{(\epsilon)}(\g)\h_{\lambda+\rho}$ the ideal of $U_{(\epsilon)}(\g)$ generated by elements of the form $X+(\lambda+\rho)(X),\;X\in\h$. Further results relate the specialization algebra  $H^0_{(\epsilon=1)}(\h_{\lambda}^{\bot},d^{(\epsilon=1)}_{\h_{\lambda}^{\bot},\mathfrak{q}})$ and other deformations of $H^0(\h_{\lambda}^{\bot},d_{\h_{\lambda}^{\bot},\mathfrak{q}})$ and $(U(\g)/U(\g)\h_{\lambda+\rho})^{\h}$.
\end{abstract}


\section{Introduction}

Our motivation arises from the last paragraphs of M. Kontsevich's paper \cite{K} which provided among several other things, a new approach to the Duflo isomorphism. The techniques of \cite{K} for the deformation quantization of a Poisson manifold $X$ were later expanded by A.S Cattaneo and G.Felder in \cite{CF2}-\cite{CF3} considering a coisotropic submanifold $C\subset X$. The theory involving two coisotropic submanifolds $C_1,C_2 \subseteq X$  with $C_1\neq C_2$, is usually termed as \textsl{biquantization}. Later, A.S. Cattaneo and Ch. Torossian in \cite{CT} set the foundations of this  perspective in the case of a Lie algebra $\g$ (and especially for the case of a symmetric space structure $\g=\mathfrak{p}\oplus\h$), considering $X=C_1=\g^{\ast}$. This paper concerns biquantization in the Lie case in general. It is meant to provide results towards a proof of the relative Duflo conjecture stated below in a way similar to the approach of the Duflo theorem in \cite{K}.

The paper is structured as follows: In Section 2 we recall the Duflo Conjecture along with the necessary techniques and results from deformation quantization and biquantization. In Section 3 we describe in detail the reduction algebra and establish some of its features that will help our arguments later. Some preliminary results used in the proof of Theorem \ref{main} and the rest of the paper are also proved in Section 4. Section 5  presents our main result, theorem \ref{main},  stating that there is a non-canonical isomorphism of associative algebras,  such that $H^0_{(\epsilon)}(\h_{\lambda}^{\bot},d^{(\epsilon)}_{\h_{\lambda}^{\bot},\mathfrak{q}})\simeq (U_{(\epsilon)}(\g)/U_{(\epsilon)}(\g)\h_{\lambda+\rho})^{\h}$.  We also describe its formula in detail. The proof is based entirely on deformation quantization techniques and the idea of translating into equations, the concentrations of configurations spaces needed to solve a Stokes equation. The idea behind the Stokes argument comes from Kontsevich's way to prove his formality theorem and associativity for his $\ast-$ product in \cite{K}. In the isomorphism formula, there is a so far unknown term which is  described in terms of colored Kontsevich graphs. This means that there is now an explicit link between the reduction algebra $H^0_{(\epsilon)}(\h_{\lambda}^{\bot},d^{(\epsilon)}_{\h_{\lambda}^{\bot},\mathfrak{q}})$, as the right candidate for the quantization of $S(\g/\h)^{\h}$ and the other side of the Duflo conjecture, namely $ (U(\g)/U(\g)\h_{\lambda+\rho})^{\h}$. Finally,  Section 6 contains results on the specialization algebra $H^0_{(\epsilon=1)}(\h_{\lambda}^{\bot},d^{(\epsilon=1)}_{\h_{\lambda}^{\bot},\mathfrak{q}}):=H^0_{(\epsilon)}(\h_{\lambda}^{\bot},d^{(\epsilon)}_{\h_{\lambda}^{\bot},\mathfrak{q}})/<\epsilon-1>$ and various deformations of $(U(\g)/U(\g)\h_{\lambda})^{\h}$. They provide an insight for deformations of $H^0(\h_{\lambda}^{\bot},d_{\h_{\lambda}^{\bot},\mathfrak{q}})$ and $(U(\g)/U(\g)\h_{\lambda+\rho})^{\h}$ and will help prove some results on characters of invariant differential operators in a subsequent paper.\newline
\noindent \textbf{Acknowledgements.} This is part of the material presented in the author's Phd thesis at Universite Paris 7. The author would like to gratefully thank Charles Torossian for his support, exciting ideas and careful supervision during these years. He would also like to thank Fred Van Oystaeyen and Simone Gutt for their kind hospitality at the universities of Antwerp and Brussels respectively.

\section{Duflo Conjecture and Deformation (bi)quantization.}

\textbf{2.1. Duflo Conjecture.} In this section we introduce some notation and state the Duflo Conjecture. A  review of the Conjecture can be found in \cite{TOR}, $\mathcal{x}$ 1.1-1.2. Let $G$ be a connected and simply connected Lie group with Lie algebra $\g$ of finite dimension. Let $\h\subset\g$ be a subalgebra with Lie group $H$, $\lambda$ a character of $\h$ and $\chi_{\lambda}:\;H\longrightarrow \mathbb{C}$ the unitary character defined by the formula $\chi_{\lambda}(\exp Y)= \exp(i\lambda(Y))$, for $Y\in\h$. Let $U(\g)$ be the universal enveloping algebra of $\g$ and $U_{\mathbb{C}}(\g):=U(\g)\otimes\mathbb{C}$ be its complexification. Denote as $U_{\mathbb{C}}(\g)\h_{i\lambda}$ the ideal of $U_{\mathbb{C}}(\g)$ generated by elements of the form $\{Y +i\lambda(Y)/ Y\in\h\}$ and set $U_{\mathbb{C}}(\g,\h,\lambda):= \{A \in U_{\mathbb{C}}(\g) / \; \forall Y\in \h,\;[A,Y] \in U_{\mathbb{C}}(\g)\h_{i\lambda}\}$. We denote as $C^{\infty}(G,H,\chi_\lambda)$ the vector space of complex smooth functions~$\theta$ on $G$ that satisfy $\theta (gh)=\chi^{-1}_{\lambda}(h)\theta(g)$, $\forall h\in H, \forall g\in G$. Finally let $\mathbb{D}(\g,\h,\lambda)$ be the algebra of linear differential operators, that leave the space $C^{\infty}(G,H,\chi_\lambda)$ invariant and commute with the left translation on $G$: $\forall g\in G, \forall D\in \mathbb{D}(\g,\h,\lambda), \forall\theta\in C^{\infty}(G,H,\chi_\lambda)$, it is $D(C^{\infty}(G,H,\chi_\lambda))\subset C^{\infty}(G,H,\chi_\lambda),\;\textlatin{and}\;D(L(g)\theta)=L(g)(D(\theta))$. Koornwinder in \cite{koorn} proved that there is an algebra isomorphism
\begin{equation}\label{Koornwider}
U(\g_{\mathbb{C}},\h,\lambda)/U_{\mathbb{C}}(\g)\h_{i\lambda}\stackrel{\sim}{\rightarrow} \mathbb{D}(\g,\h,\lambda).
\end{equation}
We will change the notation to $(U_{\mathbb{C}}(\g)/U_{\mathbb{C}}(\g)\h_{\lambda})^{\h}$ instead of $U_{\mathbb{C}}(\g,\h,\lambda)/U_{\mathbb{C}}(\g)\h_{i\lambda}$ for the algebra at the left of (\ref{Koornwider}), where the exponent means the invariants with respect to the extention of the $\mathrm{ad}\h-$ action on $U(\g)$. Let $S(\g)$ be the symmetric algebra of $\g$. We also consider it as the algebra $\mathbb{R}[\g^\ast]$ of polynomials on the dual algebra $\g^\ast$ of $\g$. Let $S_{\mathbb{C}}(\g)$ be its complexification, $\hat{\lambda}$ an extention of $\lambda$ on $\g^{\ast}$, and $S_{\mathbb{C}}(\g)\h_{\lambda}$ the ideal of $S_{\mathbb{C}}(\g)$ generated by elements of the form $\{Y +\lambda(Y)/ Y\in\h\}$. The algebra $S_{\mathbb{C}}(\g)$ has a natural Poisson structure defined for $X,Y\in\g$ by $\{X,Y\}:=[X,Y]$. It induces a Poisson structure on $(S_{\mathbb{C}}(\g)/S_{\mathbb{C}}(\g)\h_{\lambda})^{\h}$. Let $\h^{\bot}:=\{l\in\g^{\ast}/l(\h)=0\}$ and $\mathbb{C}[-\hat{\lambda}+\h^{\bot}]^{\h}$ be the Poisson algebra of $\h$-invariant complex polynomial functions on $-\lambda +\h^{\bot}$. Then $(S_{\mathbb{C}}(\g)/S_{\mathbb{C}}(\g)\h_{\lambda})^{\h}\simeq \mathbb{C}[-\hat{\lambda}+\h^{\bot}]^{\h}$ as algebras. We denote by $C_{poiss}[(S(\g)/S(\g)\h_{\lambda})^{\h}]$ and $C_{ass}[(U(\g)/U(\g)\h_{\lambda})^{\h}]$ the centers of the corresponding Poisson and associative structure of these algebras over $\mathbb{R}$. Finally, for $H\in\h$, let $\rho(H)=-\frac{1}{2}\mathrm{Tr}\;\mathrm{ad}_\g(H)$. In ~\cite{DUF1}, M. Duflo stated the following:
\newtheorem{con}{Conjecture}[section]
\begin{con}
With the previous notations, there is an algebra isomorphism,

\begin{equation}\label{duflo2}
C_{poiss}[(S(\g)/S(\g)\h_{\lambda})^{\h}]\simeq C_{ass}[(U(\g)/U(\g)\h_{\lambda+\rho})^{\h}].
\end{equation}
\end{con}
\noindent Before we continue with details on deformation quantization, we comment on the Poisson structure of $(S(\g)/S(\g)\h_{\lambda})^\h$ and the deformation $(S(\g)/S(\g)\h_{t\lambda})^\h$, and fix some more notation. Consider the algebras $(S_{\mathbb{C}}(\g)/S_{\mathbb{C}}(\g)\h_{\lambda})^{\h}$,  and $(S_{\mathbb{C}}(\g)/S_{\mathbb{C}}(\g)\h_{t\lambda})^{\h}$, $t\in\mathbb{C}^{\ast}$. If we equip the second with the Poisson bracket $t\{\cdot,\cdot\}$, where $\{\cdot,\cdot\}$ is the standard Poisson bracket on $(S_{\mathbb{C}}(\g)/S_{\mathbb{C}}(\g)\h_{t\lambda})^{\h}$, then the map $\mathrm{I}_t: P(\cdot)\mapsto P(\frac{1}{t}\cdot)$ is an isomorphism of Poisson algebras. Indeed, let $\mathcal{P}ol(-\lambda+\h^{\bot})^H$ denote the $H-$ invariant polynomial functions on $-\lambda+\h^{\bot}$. We identify this algebra with $(S_{\mathbb{C}}(\g)/S_{\mathbb{C}}(\g)\h_{\lambda})^{\h}$. Let $\mathcal{P}ol\left(t(-\lambda+\h^{\bot})\right)^H$ denote the $H-$ invariant polynomial functions on $t(-\lambda+\h^{\bot})$. We identify this algebra with $(S_{\mathbb{C}}(\g)/S_{\mathbb{C}}(\g)\h_{t\lambda})^{\h}$. Consider the map $\mathrm{I}_t:\; S_{\mathbb{C}}(\g)\longrightarrow S_{\mathbb{C}}(\g)$, $t\in \mathbb{C}^{\ast}$, defined for $X\in \g$ by $\mathrm{I}_t(X)=\frac{X}{t}$. For $H\in\h$, we have $\mathrm{I}_t(H+\lambda(H))=\frac{1}{t}(H+t\lambda(H))$ and so $S(\g)\h_{\lambda}$ is mapped onto $S(\g)\h_{t\lambda}$. Thus $(S_{\mathbb{C}}(\g)/S_{\mathbb{C}}(\g)\h_{\lambda})^{\h}$ is mapped on $(S_{\mathbb{C}}(\g)/S_{\mathbb{C}}(\g)\h_{t\lambda})^{\h}$ and $(S_{\mathbb{C}}(\g)/S_{\mathbb{C}}(\g)\h_{\lambda})^{\h}\stackrel{alg}{\simeq} (S_{\mathbb{C}}(\g)/S_{\mathbb{C}}(\g)\h_{t\lambda})^{\h}$. We denote this algebra isomorphism again as $\mathrm{I}_t$. For $P\in S_{\mathbb{C}}(\g)$, we set $P_t:=\mathrm{I}_t(P)$.  Define a Poisson structure $\{\cdot,\cdot\}_{(t)},\;t\in\mathbb{C}^{\ast}$ on $S_{\mathbb{C}}(\g)$ by $\{\cdot,\cdot\}_{(t)}:=t\{\cdot,\cdot\}$, where $\{\cdot,\cdot\}$ is the standard Poisson structure on $S_{\mathbb{C}}(\g)$. To show that $\mathrm{I}_t$ is a morphism of Poisson algebras, it suffices to verify this property on the homogeneous elements. Let $P,Q\in S_{\mathbb{C}}(\g)$ be two homogeneous elements with $\mathrm{deg}(P)=p,\mathrm{deg}(Q)=q$. Then $\mathrm{I}_t(\{P,Q\})=\frac{1}{t^{p+q-1}}\{P,Q\}=t\{P_t,Q_t\}=\{P_t,Q_t\}_{(t)}$ and thus $\mathrm{I}_t:\;\left(S_{\mathbb{C}}(\g),\{\cdot,\cdot\}\right) \longrightarrow \left(S_{\mathbb{C}}(\g),\{\cdot,\cdot\}_{(t)}\right)$ is a morphism of Poisson algebras. The Poisson structure $\{\cdot,\cdot\}_{(t)}$ on $S_{\mathbb{C}}(\g)$ induces a Poisson structure $\{\cdot,\cdot\}_{(t)}$ on $(S_{\mathbb{C}}(\g)/S_{\mathbb{C}}(\g)\h_{t\lambda})^{\h}$. Since $\mathrm{I}_t$ maps $S(\g)\h_{\lambda}$ onto $S(\g)\h_{t\lambda}$, the morphism of Poisson algebras $\mathrm{I}_t$ induces another morphism of Poisson algebras (which we denote with the same symbol) $\mathrm{I}_t:\;\left((S_{\mathbb{C}}(\g)/S_{\mathbb{C}}(\g)\h_{\lambda})^{\h},\;\{\cdot,\cdot\}\right)\longrightarrow \left((S_{\mathbb{C}}(\g)/S_{\mathbb{C}}(\g)\h_{t\lambda})^{\h},\;\{\cdot,\cdot\}_{(t)}\right)$. It is injective and surjective and thus an isomorphism of Poisson algebras. \newline
\noindent Let $\mathbf{K}=\mathbb{R}$ or $\mathbb{C}$ and $T(\g)$ be the tensor algebra of $\g$. Set $T_{(\epsilon)}(\g):=\mathbf{K}[\epsilon]\otimes T(\g)$ and let $\mathcal{I}_{\epsilon}$ be the two-sided ideal $<X\otimes Y-Y\otimes X-\epsilon[X,Y]>$ of $T_{(\epsilon)}(\g)$. Define the deformed universal enveloping algebra of $\g$ as $U_{(\epsilon)}(\g):=T_{(\epsilon)}(\g)/\mathcal{I}_{\epsilon}$. Otherwise, $U_{(\epsilon)}(\g)$ can be defined considering the Lie algebra $\g_{\epsilon}$ over $\mathbf{K}[\epsilon]$ with Lie bracket for $X,Y\in \g$ defined as $[X,Y]_{\epsilon}:=\epsilon[X,Y]$. Then $U_{(\epsilon)}(\g)$ is the universal enveloping algebra $U(\g_{\epsilon})$ over the ring $\mathbf{K}[\epsilon]$. We define similarly $S_{(\epsilon)}(\g):=T_{(\epsilon)}(\g)/<X\otimes Y- Y\otimes X>$. If $\mathcal{A}$ is an algebra, we denote the polynomials and series in $\epsilon$ with coefficients in $\mathcal{A}$ by $\mathcal{A}[\epsilon]$ and $\mathcal{A}[[\epsilon]]$ respectively.

\noindent\textbf{2.2. Deformation Quantization.} We briefly recall facts from \cite{K},\cite{CF2},\cite{CF3},\cite{CT}. Let $X$ be a Poisson manifold with $\dim(X)=n$, and $\{x_1,\ldots,x_n\}$ a system of local coordinates on $X$. Let $C^{\infty}(X)$ be the algebra of smooth functions on $X$. For a multiindex $R=(r_1,\ldots,r_p),1\leq r_i\leq n,\;1\leq i\leq p$ with $r_i\in\mathbb{N}\cup \{0\}$, let $\partial_R(f):=\frac{\partial^p f}{\partial x_{r_1}\cdots\partial x_{r_p}}$. Let also $(B_i)_{i\geq 1}$ be bidifferential operators on $C^{\infty}(X)$ meaning that the $B_i$ are  locally of bounded order and can be written in the form $B_i(f,g)=\sum_{R,S}b_i^{RS}\partial_R(f)\partial_S(g)$ with $f,g\in C^\infty(X)$. The terms $b_i^{RS}$ are smooth locally defined functions and $b_i^{RS}\neq~ 0$ for finite number of multiindices $R,S$. We suppose that $\forall i$, $B_i$ is of finite order. A star-product on $C^{\infty}(X)$ is an $\mathbb{R}[[\epsilon]]$-bilinear map $C^{\infty}(X)[[\epsilon]]\times C^{\infty}(X)[[\epsilon]]\longrightarrow C^{\infty}(X)[[\epsilon]],\; (f,g)\mapsto f\ast g$, which for $f,g,h \in C^{\infty}(X)$, satisfies:
\begin{enumerate}
\item $f\ast g= f\cdot g + \sum_{i=1}^{\infty}B_i(f,g)\epsilon^i$
\item $(f\ast g)\ast h = f\ast(g\ast h)$
\item $f\ast 1 = 1\ast f =f$.
\end{enumerate}
The most common example of a $\ast-$ product is the Moyal product: Let $(x_1,\cdots, x_k)$ be the local coordinates associated to the basis $\{e_1,\ldots,e_k\}$ of $\mathbb{R}^k$ and denote as $\partial_s$ the partial derivative with respect to the $s^{th}-$ coordinate. A Poisson structure on $\mathbb{R}^k$ is an antisymmetric $k\times k$ matrix $(\pi^{ij})\in M(\mathbb{R}^{k^2})$ and the Poisson bivector $\pi$ associated to $(\pi^{ij})$ is defined for $ f,g\in C^{\infty}(\mathbb{R}^k)$, by $\pi(f,g):=\frac{1}{2}\sum_{i,j=1}^k \pi^{ij}\partial_i(f) \otimes \partial_j(g)$. The Moyal product on $C^{\infty}(\mathbb{R}^k)$ is the star product defined by the formula
\begin{equation}\label{moyal product}
f\ast_Mg=\sum_{n=0}^{\infty}\frac{\epsilon^n}{n!}\sum_{i_1,\ldots,i_n;j_1,\ldots,j_n}\prod_{s=1}^n\pi^{i_sj_s}\left(\prod_{s=1}^n\partial_{i_s}\right)(f)\cdot\left(\prod_{s=1}^n\partial_{j_s}\right)(g)
\end{equation}
where the product $\cdot$ at the right hand side is the pointwise product on $C^{\infty}(\mathbb{R}^k)$. In the sequence we will heavily use the following graphs introduced by Kontsevich in \cite{K}: Denote by $\mathbf{Q_{n,\overline{m}}}$ the set of all \textsl{admissible} graphs $\Gamma$, that is graphs satisfying the following properties:
\begin{enumerate}
\item The set $V(\Gamma)$ of vertices of $\Gamma$ is the disjoint union of two ordered sets $V_1(\Gamma)$ and $V_2(\Gamma)$, isomorphic to $\{1,\ldots,n\}$ and $\{1,\ldots,\overline{m}\}$ respectively. Their elements are called \textsl{type I} vertices, for $V_1(\Gamma)$, and \textsl{type II} vertices, for $V_2(\Gamma)$. 
\item The number of type I and type II vertices must satisfy the  inequalities $n,\overline{m}\geq 0,\;\; 2n+\overline{m}-2\geq 0$.
\item The set $E(\Gamma)$ of edges of the graph is finite. Each edge starts from a type I vertex and ends to a vertex of type I or type II. No loops or double edges are allowed for $\Gamma$.
\item All elements of $E(\Gamma)$ are oriented and the set of edges $S(r)$ starting from $r\in V_1(\Gamma)$ is ordered.
\item The set $E(\Gamma)$ is ordered in a compatible way with the order in $V_1(\Gamma)$, and $S(r)$ for $r\in V_1(\Gamma)$.
\end{enumerate}
To each such graph one associates a bidifferential operator: Let $\psi_1,\ldots,\psi_n$ be $n$ multivector fields on $\mathbb{R}^k$ and suppose that each $\psi_r$ is a skew symmetric tensor of degree $k_r$. To the tensor product $\psi_1\otimes\ldots\otimes\psi_n$ and to a graph $\Gamma\in \mathbf{Q_{n,\overline{m}}}$ we associate a differential operator $B_{\Gamma}$ as follows: Use the notation $[[1,k]]:=\{1,\ldots,k\}$ and let $L:\;E(\Gamma)\longrightarrow [[1,k]]$ be a labelling function.
\begin{enumerate}
\item Fix a vertex $r\in [[1,n]]$. If $card(S(r))\neq k_r$, set $B_{\Gamma}=0$. If $card(S(r))= k_r$, let $S(r)=\{e_r^1,\ldots,e_r^{k_r}\}$ be the ordered set of edges leaving $r$. Associate the function $\psi_r^{L(e_r^1),\ldots,L(e_r^{k_r})}$ to $r$ .
\item On each vertex $\overline{1},\ldots,\overline{m}\in V_2(\Gamma)$ we associate respectively a function $F_1,\ldots, F_m\in C^{\infty}(\mathbb{R}^k)$. 
\item  To the $p^{th}-$ edge of $S(r)$, associate the partial derivative w.r.t the coordinate variable $L(e_r^p)$.
\item This derivative acts on the  function associated to $v\in V_1(\Gamma)\cup V_2(\Gamma)$ where the edge $e_r^p$ arrives.
\end{enumerate}
Since $E(\Gamma)\subset V_1(\Gamma)\times (V_1(\Gamma)\cup V_2(\Gamma))$, let $(p,m)\in E(\Gamma)$ represent an oriented edge of $\Gamma$ from $p$ to $m$. In this paper we deal with at most two type II vertices, so we restrict to the case $card(V_2(\Gamma))=2$ and set hereafter  $F,G$ to be two functions corresponding to the ordered set $V_2(\Gamma)$.  In the Lie case it is $\psi_i=\pi=[\cdot,\cdot]$, so we need $card(S(r))=2$, $\forall  r\in V_1(\Gamma)$. The operator associated to a $\Gamma\in\mathbf{Q_{n,2}}$ and $\pi\otimes\cdots\otimes\pi$ (n-times) is defined as
\begin{equation}\label{bidiff op}
B_{\Gamma}^{\pi}(F,G)= \sum_{L:E(\Gamma)\rightarrow [[1,k]]}\left[\prod_{r=1}^{\#(V_1(\Gamma))}\left(\prod_{\delta\in E(\Gamma),\;\delta=(\cdot,r)}\partial_{L(\delta)}\right)\pi^{L(e_r^1)L(e_r^2)}\right]\times
\end{equation}
\[\times\left(\prod_{\delta\in E(\Gamma),\;\delta=(\cdot,\overline{1})}\partial_{L(\delta)}\right)(F)\times\left(\prod_{\delta\in E(\Gamma),\;\delta=(\cdot,\overline{2})}\partial_{L(\delta)}\right)(G).\]
Finally we recall the last ingredient for Kontsevich's $\ast-$ product formula, the coefficient $\omega_{\Gamma}$. Let $\mathcal{H}=\{z\in \mathbb{C}/ \mathfrak{Im}(z)\geq 0\}$ be the upper-half plane and let $\mathcal{H}^+=\{z\in \mathbb{C}/ \mathfrak{Im}(z)> 0\}$. Embed an admissible graph $\Gamma$ in $\mathcal{H}$ by putting the type II vertices on the real axis (they can move on the axis) and letting the type I vertices move in $\mathcal{H}^+$. Consider the configuration space $C_{n,\overline{m}}$ defined as
\[C_{n,\overline{m}}:=\{(z_1,\ldots,z_n,z_{\overline{1}},\ldots,z_{\overline{m}})\in \mathbb{C}^{n+\overline{m}}/z_i\in \mathcal{H}^+,\;\;z_{\overline{\i}} \in \mathbb{R},\;z_i\neq z_j \;\textlatin{for}\; i\neq j,\;z_{\overline{\i}}\neq z_{\overline{\j}} \;\textlatin{for} \;\overline{\i}\neq \overline{\j}\}.\]
For our purposes we need the quotient $\overline{C}_{n,\overline{m}}:=C_{n,\overline{m}}/G_2$, where $G_2$ is the 2-dimensional Lie group of horizontal translations and dilations in $\mathcal{H}$. The action is $G_2 \times C_{n,\overline{m}}\longrightarrow C_{n,\overline{m}}$, $<(a,b)\times z>\mapsto az+b$ where $a\in\mathbb{R}^+,b\in\mathbb{R}$. It is a free action so since $C_{n,\overline{m}}$ is a smooth manifold, so will be $\overline{C}_{n,\overline{m}}$ with $\dim_{\mathbb{R}}(\overline{C}_{n,\overline{m}})=2n+\overline{m}-2$. The final restriction is made so that the configuration manifold we deal with is connected. For this let $\overline{C}_{n,\overline{m}}^+$ be the connected component of $\overline{C}_{n,\overline{m}}$ obtained by putting an ordering on the type II vertices. In other words, $\overline{C}_{n,\overline{m}}^+:=\{(z_1,\ldots,z_n,z_{\overline{1}},\ldots,z_{\overline{m}})\in\overline{C}_{n,\overline{m}}/\forall i<j,\;\; z_{\overline{\i}}<z_{\overline{\j}}\}$. Consider now the manifold $\overline{C}_{2,0}$ and a map on it, called the \textsl{angle map}, defined by $\phi(z_1,z_2):=\frac{1}{2\pi}\mathrm{arg}(\frac{z_1-z_2}{\overline{z}_1-z_2}),\;z_1\neq z_2$, where $\mathrm{arg}(\bullet)$ denotes an argument function. This map descends to $\overline{C}_{2,0}$. Let then $e=(z_i,z_j)$ be an edge of $\Gamma$, and consider the natural projection $p_{e}: \overline{C}_{n,\overline{m}}\longrightarrow \overline{C}_{2,0},\;(z_1,\ldots,z_n,z_{\overline{1}},\ldots,z_{\overline{m}})\mapsto (z_i,z_j)$. The pullback $p^{\ast}_{e}$ applied on the 1-form $\mathrm{d}\phi$ defines a form on $\overline{C}_{n,\overline{m}}$, namely $p^{\ast}_{e}(\mathrm{d}\phi)=:\mathrm{d}\phi_{e} \in \Omega^1(\overline{C}_{n,\overline{m}})$. Then define $\Omega_{\Gamma}$ to be the form $\Omega_{\Gamma}:=\bigwedge_{e\in E(\Gamma)}\mathrm{d}\phi_e$. In fact there is an order with which the 1-forms $\mathrm{d}\phi_e$ appear in this exterior product. This order is the one induced by the graph $\Gamma$, considering first the order in $V_1(\Gamma)$ and then the order in $S(r),\; r\in V_1(\Gamma)$. To each $\Gamma \in \mathbf{Q_{n,\overline{m}}}$ Kontsevich associates a coefficient $\omega_{\Gamma}\in\mathbb{R}$ by the formula $\omega_{\Gamma}:=\frac{1}{(2\pi)^{2n+\overline{m}-2}}\int_{\overline{C}_{n,\overline{m}}^+}\Omega_{\Gamma}$. 
Since $\dim(\overline{C}_{n,\overline{m}}^+)=2n+\overline{m}-2$, this is well defined iff $\Omega_{\Gamma}$ is a ($2n+\overline{m}-2$)-form. As for convergence, Kontsevich has constructed a compactification of $\overline{C}_{n,\overline{m}}^+$ to which $\Omega_{\Gamma}$ extends continuously. We consider graphs $\Gamma$ such that a priori $\omega_{\Gamma}\neq 0$. 
\newtheorem{ok}[con]{Theorem}
\begin{ok} \cite{K} \label{kon}
Let $\pi$ be a Poisson bivector on $\mathbb{R}^k$ such that $(\mathbb{R}^k,\pi)$ is a Poisson manifold. Then for $f,g\in C^{\infty}(\mathbb{R}^k)$, the operator $f\ast_Kg:=fg+\sum_{n=1}^{\infty}\epsilon^n\left(\frac{1}{n!}\sum_{\Gamma\in\mathbf{Q_{n,2}}}\omega_{\Gamma}B_{\Gamma}^{\pi}(f,g)\right)$ is an associative product. 
\end{ok}
This is the local result while the global result for a Poisson manifold $X$ was given in \cite{CFT}.

\noindent\textbf{2.3. Coisotropic submanifolds.}  Let $(X,\pi)$ be a Poisson manifold, and $\mathcal{T}X$, $\mathcal{T}^\ast X$ the tangent and cotangent bundles respectively. Denote by $\mathrm{d}g\in\mathcal{T}^{\ast}X$ the covector field associated to a function $g\in C^{\infty}(X)$ by the relation $\mathrm{d}g(L)(s)=(L_s(g))(s)$ for $s\in X$, $L\in\mathcal{T}X$ and $L_s$ its value on $s$. Let now $C\subset X$ be a submanifold of $X$, $\mathcal{T}C$ the tangent bundle of $C$, and $N^{\ast}C= Ann(\mathcal{T}C)\subset\mathcal{T}^{\ast}X|_{C}$ the conormal bundle of $C$. Let $\pi^{\#}:\mathcal{T}^{\ast}X\longrightarrow \mathcal{T}X,\;\;Z\mapsto \pi^{\#}(Z)$ defined for $Y\in \mathcal{T}^{\ast}X,s\in C$ by the formula $<\pi^{\#}_s(Z_s),Y_s>=<\pi_s,Z_s\otimes Y_s>$, be the bundle map on the cotangent space induced by $\pi$. The submanifold $C$ is called \textsl{coisotropic}, if $\pi^{\#}|_{C}(N^{\ast}C)\subset \mathcal{T}C$. In other words $C$ is coisotropic if the ideal $I(C)\subset C^{\infty}(X)$ of functions vanishing on $C$ is a Poisson subalgebra of $C^{\infty}(X)$. It is clear that in the Lie case, for a subalgebra $\h\subset\g$, the annihilator $\h^{\bot}$ (and also $\g^{\ast}$) is a coisotropic submanifold of $\g^{\ast}$ with the natural Poisson structure induced by the Lie bracket. In two fundamental papers ~\cite{CF2},~\cite{CF3}, A. Cattaneo and G. Felder, expanded the formality theorem of M. Kontsevich for the case of one coisotropic submanifold $C$ of a Poisson manifold $X$. More specifically, the main theorem of ~\cite{CF3} states that the DGLA of multivector fields on an infinitesimal neighborhood of $C$ is $L_{\infty}-$quasiisomorphic to the DGLA of multidifferential operators on $\Gamma(C,\wedge NC)$, the sections of the exterior algebra of the conormal bundle of $C$. Their idea was to replace the Poisson algebra $C^{\infty}(X)$ with the graded commutative algebra $\mathcal{A}=\Gamma(C,\wedge NC)$ and then use a Fourier transform to compose with the original $L_{\infty}-$quasi-isomorphism constructed by Kontsevich. 

\noindent We need some facts from the local construction described in \cite{CF2}. primarily interested in the notions of \textsl{colors}, \textsl{reduction algebra} and a certain bimodule structure; we first recall the notion of colors.  In deformation quantization of a (linear) Poisson manifold, one trivially has a single \textlatin{color} for every edge in a graph $\Gamma$ since for  $e\in E(\Gamma)$, $L(e)$ determines a basis variable of $\g$ without any discrimination. In the case of (trivial) biquantization where $C_1=\g^\ast$ and $C_2=\h^\bot$ however we consider two colors, with respect to $\h$. Now each edge of a colored graph $\Gamma$ carries a color, either $(+)$ or $(-)$. Double edges are not allowed, meaning edges with the same color, source and target. We make that explicit: Let $\g$ be a Lie algebra, consider $\g^{\ast}$ as a Poisson manifold. It is easy to see that $C^{\infty}(\g^{\ast})=S(\g)$, since the symmetric algebra $S(\g)$ can be regarded as the polynomial functions on $\g^{\ast}$. Let $\mathfrak{q}$ be a supplementary space of $\h$, that is $\g=\h\oplus\mathfrak{q}$ and let $\{H_1,H_2,\ldots, H_t\}$ be a basis for $\h$ and $\{Q_1,\ldots, Q_r\}$ a basis for $\mathfrak{q}$ forming together a basis for $\g$. We identify spaces $\mathfrak{q}^{\ast}\simeq\g^{\ast}/\h^{\ast}\simeq\h^{\bot}$.  For $e\in E(\Gamma)$, let $c_e\in\{+,-\}$ be its color. Let $L:\;E(\Gamma)\longrightarrow \{1,\ldots,t,t+1,\ldots t+r\}$, satisfying $L(e)\in\{1,\ldots,t\}\;\;\textlatin{if}\;c_e=-\;\;,\;\;L(e)\in\{t+1,\ldots,t+r\}\;\;\textlatin{if}\;c_e=+$ be a 2-colored labelling function. This way, the dual basis variables $\{H_1^{\ast},H_2^{\ast},\ldots, H_t^{\ast}\}$ of $\h^{\ast}$ are associated to the color $(-)$ and dual basis variables $\{Q_1^{\ast},\ldots, Q_r^{\ast}\}$ of $\mathfrak{q}^{\ast}$ are associated to $(+)$. Graphically, the color $(-)$ will be represented with a dotted edge and the color $(+)$ will be represented with a straight edge. Recall that in order to construct the differential operator $B_{\Gamma}$ corresponding to a given graph, we associated a coordinate variable to each edge of the graph. 
In terms of $\g^{\ast}$, let $(x_i^{\ast})_{i=1,\ldots,n}$ be the coordinates relatively to the basis $\{H_1^{\ast},\ldots,H_t^{\ast},Q_1^{\ast},\ldots,Q_r^{\ast}\}$ and let $\Gamma$ be an admissible graph with two colors. The formula (\ref{bidiff op}) of $B_{\Gamma}^\pi$ in this case has to be modified in the sense that for $F,G\in S(\mathfrak{q})$, the correct formula is the same as in (\ref{bidiff op}) but using the 2-colored labelling function $L$ that we just described.

\noindent The computation of Kontsevich's coefficients is also modified: To every 2-colored graph $\Gamma$ is associated an 1-form $\overline{\Omega}_{\Gamma}$ and a coefficient $\overline{\omega}_{\Gamma}$ (the bar is used to indicate the existence of colors) as follows. Set $\phi_+,\phi_- $ to be the functions
$\phi_+(z_1,z_2):=\phi(z_1,z_2)$ and $\phi_-:=\phi(z_2,z_1)$, the bar standing here for the complex conjugate, or alternatively, $\phi_+(z_1,z_2)=\mathrm{arg}(z_1-z_2)+\mathrm{arg}(z_1-\overline{z}_2),\;\phi_{-}(z_1,z_2)=\mathrm{arg}(z_1-z_2)-\mathrm{arg}(z_1-\overline{z}_2)$. These functions also descend to $\overline{C}_{2,0}=C_{2,0}/G_2$ and thus can be used for definitions analogous to Kontsevich's. The function $\phi_+$ will be used when the edge from $z_1$ to $z_2$ carries a tangent variable (in $\mathfrak{q}^{\ast}$) and the function $\phi_-$ when the edge carries a normal variable (in $\h^{\ast}$). The form $\overline{\Omega}_{\Gamma}$ of a 2-colored graph $\Gamma$ is similarly defined as $\overline{\Omega}_{\Gamma}:=\wedge_{e\in E(\Gamma)}\mathrm{d}\phi_{\cdot,e}$ where $\mathrm{d}\phi_{+,e}=p_e^{\ast}(\mathrm{d}\phi_+),\;\mathrm{d}\phi_{-,e}=p_e^{\ast}(\mathrm{d}\phi_-)$, when $c_e=(+)$ or $(-)$ respectively. The colored coefficient is $\overline{\omega}_{\Gamma}:=\frac{1}{(2\pi)^{2n+\overline{m}-2}}\int_{\overline{C}_{n,\overline{m}}^+}\overline{\Omega}_{\Gamma}$. Let $\mathbf{Q}_{n,2}^{(2)}$ denote the set of admissible graphs with two colors and two type II vertices. Then theorem \ref{kon} is generalized in the following sense (\cite{CF3}): Consider $\mathbb{R}^r\subset \mathbb{R}^{t+r}$ as a coisotropic submanifold of $\mathbb{R}^{t+r}$. Cattaneo and Felder associate a curved $A_{\infty}$ algebra, which in the linear Poisson case is flat. Its $0-$th cohomology then accepts an associative product on it, called the \textsl{Cattaneo-Felder} product $\ast_{CF,\epsilon}:\;C^{\infty}(\mathbb{R}^r)\times C^{\infty}(\mathbb{R}^r)\longrightarrow C^{\infty}(\mathbb{R}^r)$. It is given by the formula $F\ast_{CF,\epsilon} G:=F\cdot G+\sum_{n=1}^{\infty}\epsilon^n\left(\frac{1}{n!}\sum_{\Gamma\in\mathbf{Q^{(2)}_{n,2}}}\overline{\omega}_{\Gamma}B_{\Gamma}^{\pi}(F,G)\right)$.

\section{Reduction algebras.}
\textbf{3.1 Graphs and differentials.} We specify some particular colored graphs that we will use (see ~\cite{CT} $\mathcal{x}$ 1.3, 1.6 and ~\cite{BAT} $\mathcal{x}$ 2.3). They are colored graphs with an edge colored by $(-)$ which has no end and that $\#V_2(\Gamma)=1$, i.e there is only one type II vertex. We will say that this edge "points to $\infty$" and denote it by $e_{\infty}$.  The set of such graphs with $n$ type I vertices will be denoted by $\mathbf{Q}^{\infty}_{n,1}$.
\newtheorem{de}[con]{Definition}
\begin{de}
\begin{enumerate}
\item \textbf{Bernoulli.} The Bernoulli type graphs with $i$ type I vertices, $i\in\mathbb{N}$, will be denoted by $\mathcal{B}_i$. They derive the function $F$ $i$ times, have $2i$ edges and leave an edge towards $\infty$. These conditions imply the existence of a vertex $s\in V_1(\Gamma)$ that receives no edge, called the \textsl{root} of $\Gamma$.
\item \textbf{Wheels.} The wheel type graphs with $i$ type I vertices, $i\in\mathbb{N}$, will be denoted by $\mathcal{W}_i$. They derive the function $F$ $i$ times, have $2i$ edges and leave no edge to $\infty$. 
\item \textbf{Bernoulli attached to a wheel.} Graphs of this type with $i$ type I vertices, $i\in\mathbb{N}$, will be denoted by $\mathcal{BW}_i$. They derive the function $F$ $i-1$ times and leave an edge to $\infty$.
For an $\mathcal{W}_m-$ type graph $W_m$ attached to a $\mathcal{B}_l-$ type graph $B_l$, we will write $B_lW_m\in\mathcal{B}_l\mathcal{W}_m$. Obviously $\mathcal{B}_l\mathcal{W}_m\subset\mathcal{BW}_{l+m}$.
\end{enumerate}
\end{de}
Let us now give the definition of the reduction algebra without character. Let $\{e^1_l,e^2_l\}$ be the ordered set of edges leaving the vertex $l\in V_1(\Gamma)$ of a colored graph $\Gamma \in\mathbf{Q}^{\infty}_{s,1}$. For such a $\Gamma$ and using the notation $H_i^{\ast}:=\partial_i$, let $B_{\Gamma}:\;S(\mathfrak{q})\longrightarrow S(\mathfrak{q})\otimes \h^{\ast},\;F\mapsto B_{\Gamma}(F)$ be the operator defined by the formula
\begin{equation}
B_{\Gamma}(F)= \sum_{\substack{L:E(\Gamma)\rightarrow [[1,t+r]]\\ \textlatin{L colored}\\}}\left[\prod_{r=1}^{n}\left(\prod_{e\in E(\Gamma),\;e=(\cdot,r)}\partial_{L(e)}\right)\pi^{L(e_r^1)L(e_r^2)}\right]\times\left( \prod_{\substack{e\in E(\Gamma)\\ e=(\cdot,\overline{1})\\}}\partial_{L(e)}F\right)\otimes H^{\ast}_{L(e_{\infty})}
\end{equation}
Equivalently, one writes $B_{\Gamma}(F)=\sum_{i=1}^t B_i(F) \cdot H_i^{\ast}\in S(\mathfrak{q})\otimes \h^{\ast}$, where
 \begin{equation}
 B_i(F)= \sum_{\substack{L:E(\Gamma)\rightarrow [[1,t+r]]\\ \textlatin{L colored}\\L(e_{\infty})=i\\}}\left[\prod_{r=1}^{n}\left(\prod_{\in E(\Gamma),\;e=(\cdot,r)}\partial_{L(e)}\right)\pi^{L(e_r^1) L(e_r^2)}\right]\times \left( \prod_{\substack{e\in E(\Gamma)\\ e=(\cdot,\overline{1})\\}}\partial_{L(e)}F\right)
 \end{equation}
 \newtheorem{demek}[con]{Definition}
\begin{demek}
 We denote as $d^{(\epsilon)}_{\h^{\bot},\mathfrak{q}}:\;S(\mathfrak{q})[\epsilon]\longrightarrow S(\mathfrak{q})[\epsilon]\otimes\h^{\ast}$ the differential operator $d^{(\epsilon)}_{\h^{\bot},\mathfrak{q}}=\sum_{i=1}^{\infty}\epsilon^i d^{(i)}_{\h^{\bot},\mathfrak{q}}$ where $d^{(i)}_{\h^{\bot},\mathfrak{q}}=\sum_{\Gamma\in \mathcal{B}_i\cup \mathcal{BW}_i}\overline{\omega}_{\Gamma}B_{\Gamma}$.
We define the reduction algebra $H^0_{(\epsilon)}(\h^{\bot},d^{(\epsilon)}_{\h^{\bot},\mathfrak{q}})$ as the vector space of solutions $F_{(\epsilon)}\in S(\mathfrak{q})[\epsilon]$ of the equation 
\begin{equation}\label{reduction equations}
d^{(\epsilon)}_{\h^{\bot},\mathfrak{q}}(F_{(\epsilon)})=0,
\end{equation}
 equipped with the $\ast_{CF,\epsilon}-$ product, (which is associative on $H^0_{(\epsilon)}(\h^{\bot},d^{(\epsilon)}_{\h^{\bot},\mathfrak{q}})$ by ~\cite{CF3}).
\end{demek}
\begin{figure}[h!]
\begin{center}
\includegraphics[width=9cm]{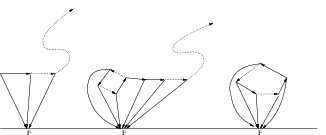}
\caption\footnotesize{A $\mathcal{B}_3$-type graph in $d^{(3)}_{\h^{\bot},\mathfrak{q}}$, a $\mathcal{B}_3\mathcal{W}_4$-type graph in $d^{(7)}_{\h^{\bot},\mathfrak{q}}$, and a $\mathcal{W}_5$-type graph. The first two are also examples of graphs in  $\mathbf{Q}^{\infty}_{3,1}$ and  $\mathbf{Q}^{\infty}_{7,1}$ respectively.}
\end{center}
\end{figure}
\textbf{3.2  The algebra} $\mathbf{H^0(\h^{\bot},d_{\h^{\bot},\mathfrak{q}})}$. If a polynomial $G\in S(\mathfrak{q})$ is homogeneous of degree $p$ with respect to the ordinary polynomial degree in $\mathfrak{q}$, we will write $\mathrm{deg}_{\mathfrak{q}}(G)=p$. Similarly we consider the $\epsilon-$ degree $\mathrm{deg}_{\epsilon}$ for elements of $S(\g)[\epsilon]$. For $F\in S(\mathfrak{q})[\epsilon]$ set $\mathrm{deg}_{\mathfrak{q},\epsilon}(F):=\mathrm{deg}_{\mathfrak{q}}(F)+\mathrm{deg}_{\epsilon}(F)$. We consider the corresponding notions of degree also for differential operators on $S(\mathfrak{q})[\epsilon]$. Let now $F_{(\epsilon)}=F_n+\epsilon F_{n-1}+\cdots +\epsilon^n F_0$, $F_i\in S(\mathfrak{q})$. The defining system (\ref{reduction equations}) of linear partial differential equations is $d^{(1)}_{\h^{\bot},\mathfrak{q}}(F_n)=0$, $d^{(2)}_{\h^{\bot},\mathfrak{q}}(F_n) + d^{(1)}_{\h^{\bot},\mathfrak{q}}(F_{n-1})=0$ and etc. That is, one uses the $\mathrm{deg}_{\epsilon}$ of the terms $\epsilon^id^{(i)}_{\h^{\bot},\mathfrak{q}}, \epsilon^{n-i}F_i$  to write down homogeneous equations. In fact it is possible to write such homogeneous equations without using  the $\mathrm{deg}_{\epsilon}$ of $\epsilon^id^{(i)}_{\h^{\bot},\mathfrak{q}}, \epsilon^{n-i}F_i$, but using instead the degree $\mathrm{deg}_{\mathfrak{q}}$ of $d^{(i)}_{\h^{\bot},\mathfrak{q}}, F_i$. This will produce the same system of homogeneous equations for a function $F=\sum_{i=0}^nF^{(n-i)}$ in $S(\mathfrak{q})$ with each $F^{(k)}$ being a homogeneous polynomial of $\mathrm{deg}_{\mathfrak{q}}(F^{(k)})=k$. We take the next few lines to explain why.\\
\noindent Recall that the possible graphs in $d^{(\epsilon)}_{\h^{\bot},\mathfrak{q}}$ are of $\mathcal{B}-$ type and $\mathcal{BW}-$ type. The root of a $\mathcal{B}-$ type graph $\Gamma$ can belong either to $\h^{\ast}$ or to $\mathfrak{q}^{\ast}$. In the first case, the contribution of the operator $B_{\Gamma}$ in $d^{(\epsilon)}_{\h^{\bot},\mathfrak{q}}$ is zero after evaluating everything at $\mathfrak{q}^{\ast}$, while in the second case the evaluation is trivial and non-zero. We suppose hereafter that these roots belong in $\mathfrak{q}$.  From a graph $B_i\in\mathcal{B}_i$ we write a differential operator that derives the function $F$ $i$ times. Since the root of the graph represents a bracket, it adds a variable (and thus a degree to $\mathrm{deg}_{\mathfrak{q}}$) to the symbol of the operator so we deduce that such operators are of degree $\mathrm{deg}_{\mathfrak{q}}(B_{B_i})=-i+1$.  Similarly, a $\mathcal{BW}-$ type graph $B_lW_m$ with $i=l+m$ type I vertices,  corresponds to a differential operator that derives the function $F$, $i-1$ number of times. Thus it is an operator of degree $\mathrm{deg}_{\mathfrak{q}}(B_{B_lW_m})=-i+1$. This means that in both cases ($\mathcal{B}$ and $\mathcal{BW}$ graphs), the operator $B_{\Gamma_i}$ for $\Gamma_i\in\mathcal{B}_i$ or $\Gamma_i\in\bigcup_{l+m=i}\mathcal{B}_l\mathcal{W} _m$, has a degree $\mathrm{deg}_{\mathfrak{q}}(B_{\Gamma_i})=-i+1$, $i$ being the number of its type I vertices. We can now regroup these graphs with respect to the degree of their associated operator and rewrite a system as (\ref{reduction equations}) but lacking the parameter $\epsilon$. \newline
\noindent Denote as $d_{\h^{\bot},\mathfrak{q}}:\;S(\mathfrak{q})\longrightarrow S(\mathfrak{q})\otimes \h$ the differential operator $d_{\h^{\bot},\mathfrak{q}}= \sum_{i=1}^{\infty}d^{(i)}_{\h^{\bot},\mathfrak{q}}$ where $d^{(i)}_{\h^{\bot},\mathfrak{q}}=\sum_{\Gamma\in \mathcal{B}_i\cup \mathcal{BW}_i}\overline{\omega}_{\Gamma}B_{\Gamma}$ and as $H^0(\h^{\bot},d_{\h^{\bot},\mathfrak{q}})$ the vector space of polynomials $F\in S(\mathfrak{q})$, solutions of the equation $d_{\h^{\bot},\mathfrak{q}}(F)=0$. In the linear Poisson case as here, the Cattaneo-Felder construction without $\epsilon$ is still valid for polynomial functions and defines an associative product $\ast_{CF}$ on $H^0(\h^{\bot},d_{\h^{\bot},\mathfrak{q}})$. In the sequence $H^0(\h^{\bot},d_{\h^{\bot},\mathfrak{q}})$ will stand for the algebra $\left(H^0(\h^{\bot},d_{\h^{\bot},\mathfrak{q}}),\ast_{CF}\right)$.

\newtheorem{p}[con]{Proposition}
\begin{p}\label{homogeneous}
Let $\g$ be a Lie algebra, $\h\subset\g$ a subalgebra and $\mathfrak{q}$ such that $\g=\h\oplus\mathfrak{q}$. If $F^{(i)}\in S(\mathfrak{q}),\;\mathrm{deg}_{\mathfrak{q}}(F^{(i)})=i$ and $F=\sum_{i=0}^n F^{(n-i)}\in H^0(\h^{\bot},d_{\h^{\bot},\mathfrak{q}})$, then $F_{(\epsilon)}:=\sum_{i=0}^n \epsilon^i F^{(n-i)}\in H^0_{(\epsilon)}(\h^{\bot},d^{(\epsilon)}_{\h^{\bot},\mathfrak{q}})$ and $\mathrm{deg}_{\mathfrak{q},\epsilon}(F_{(\epsilon)})=n$. Conversely, if $F_{\epsilon}=\sum_{j=0}^n \epsilon^jF_j \in   H^0_{(\epsilon)}(\h^{\bot},d^{(\epsilon)}_{\h^{\bot},\mathfrak{q}})$, then $F=\sum_{j=0}^nF_j \in H^0(\h^{\bot},d_{\h^{\bot},\mathfrak{q}})$.
\end{p}
\textit{Proof.} As explained, the algebra $H^0(\h^{\bot},d_{\h^{\bot},\mathfrak{q}})$ is defined by a system of homogeneous equations because the operators in $d_{\h^{\bot},\mathfrak{q}}$ are homogeneous with respect to $\mathrm{deg}_{\mathfrak{q}}$. So these are
\begin{equation}\label{fst}
d^{(1)}_{\h^{\bot},\mathfrak{q}}(F^{(n)})=0,\;\;d^{(2)}_{\h^{\bot},\mathfrak{q}}(F^{(n)})+d^{(1)}_{\h^{\bot},\mathfrak{q}}(F^{(n-1)})=0, \;\;d^{(3)}_{\h^{\bot},\mathfrak{q}}(F^{(n)})+d^{(2)}_{\h^{\bot},\mathfrak{q}}(F^{(n-1)})+d^{(1)}_{\h^{\bot},\mathfrak{q}}(F^{(n-2)})=0
\end{equation}
and so on. The terms $F^{(i)}$ are homogeneous polynomials of degree $\mathrm{deg}_{\mathfrak{q}}(F^{(i)})=i$ and we need to show that the equations defining $H^0_{(\epsilon)}(\h^{\bot},d^{(\epsilon)}_{\h^{\bot},\mathfrak{q}})$ are satisfied by $F_{(\epsilon)}=\sum_{i=0}^n \epsilon^i F^{(n-i)}$. Indeed they are homogeneous with respect to the total degree $\mathrm{deg}_{\mathfrak{q},\epsilon}:=\mathrm{deg}_{\mathfrak{q}}+\mathrm{deg}_{\epsilon}$. More precisely and using temporarily the notation $\epsilon^{i}d^{(i)}_{\h^{\bot},\mathfrak{q}}=: d^{\epsilon,i}_{\h^{\bot},\mathfrak{q}}$ we have by definition $d^{(\epsilon)}_{\h^{\bot},\mathfrak{q}}=\sum  d^{\epsilon,i}_{\h^{\bot},\mathfrak{q}}$. Then for $F_{(\epsilon)}=\sum_{i=0}^n \epsilon^i F^{(n-i)}$,
\[d^{(\epsilon)}_{\h^{\bot},\mathfrak{q}}(F_{(\epsilon)})=d^{\epsilon,1}_{\h^{\bot},\mathfrak{q}}(F^{(n)})+d^{\epsilon,1}_{\h^{\bot},\mathfrak{q}}(\epsilon F^{(n-1)})+d^{\epsilon,2}_{\h^{\bot},\mathfrak{q}}(F^{(n)})+d^{\epsilon,1}_{\h^{\bot},\mathfrak{q}}(\epsilon^2F^{(n-2)})+d^{\epsilon,2}_{\h^{\bot},\mathfrak{q}}(\epsilon F^{(n-1)})+d^{\epsilon,3}_{\h^{\bot},\mathfrak{q}}(F^{(n)})+\cdots=\] 
\[\epsilon\left(d^{(1)}_{\h^{\bot},\mathfrak{q}}(F^{(n)})\right)+\epsilon^2\left(d^{(2)}_{\h^{\bot},\mathfrak{q}}(F^{(n)})+d^{(1)}_{\h^{\bot},\mathfrak{q}}(F^{(n-1)})\right)+\epsilon^3\left(d^{(3)}_{\h^{\bot},\mathfrak{q}}(F^{(n)})+d^{(2)}_{\h^{\bot},\mathfrak{q}}(F^{(n-1)})+d^{(1)}_{\h^{\bot},\mathfrak{q}}(F^{(n-2)})\right)+\cdots\]
which by the system (\ref{fst}), gives $d^{(\epsilon)}_{\h^{\bot},\mathfrak{q}}(F_{(\epsilon)})=0$ and so $F_{(\epsilon)}\in H^0_{(\epsilon)}(\h^{\bot},d^{(\epsilon)}_{\h^{\bot},\mathfrak{q}})$. Conversely, if $F_{\epsilon}=\sum_{j=0}^n \epsilon^jF_j \in H^0_{(\epsilon)}(\h^{\bot},d^{(\epsilon)}_{\h^{\bot},\mathfrak{q}})$, then $d^{(\epsilon)}_{\h^{\bot},\mathfrak{q}}(F_{\epsilon})=0$ and the first equations of this system are 

\begin{equation}\label{second system}
d^{(1)}_{\h^{\bot},\mathfrak{q}}(F_0)=0,\;\;d^{(2)}_{\h^{\bot},\mathfrak{q}}(F_0)+d^{(1)}_{\h^{\bot},\mathfrak{q}}(F_{1})=0,\;\;d^{(3)}_{\h^{\bot},\mathfrak{q}}(F_0)+d^{(2)}_{\h^{\bot},\mathfrak{q}}(F_{1})+d^{(1)}_{\h^{\bot},\mathfrak{q}}(F_{2})=0
\end{equation}
and so on. Then 
\[d_{\h^{\bot},\mathfrak{q}}(F)=d^{(1)}_{\h^{\bot},\mathfrak{q}}(F_0)+\left(d^{(2)}_{\h^{\bot},\mathfrak{q}}(F_0)+d^{(1)}_{\h^{\bot},\mathfrak{q}}(F_{1})\right)+\left(d^{(3)}_{\h^{\bot},\mathfrak{q}}(F_0)+d^{(2)}_{\h^{\bot},\mathfrak{q}}(F_{1})+d^{(1)}_{\h^{\bot},\mathfrak{q}}(F_{2})\right)+\left(\cdots\right)+\cdots\]
and from the system (\ref{second system}) we see that $F=\sum_{i=0}^nF_i$ satisfies $d_{\h^{\bot},\mathfrak{q}}(F)=0$. Thus $F\in H^0(\h^{\bot},d_{\h^{\bot},\mathfrak{q}})$. $\diamond$
 
As a corollary, one gets that the algebra $(H^0_{(\epsilon)}(\h^{\bot},d^{(\epsilon)}_{\h^{\bot},\mathfrak{q}}),\ast_{CF,\epsilon})$ is graded and homogeneous with respect to $\mathrm{deg}_{\mathfrak{q},\epsilon}$. For $\mathrm{deg}_{\mathfrak{q},\epsilon}=N$, let $\left(H^0_{(\epsilon)}(\h^{\bot},d^{(\epsilon)}_{\h^{\bot},\mathfrak{q}})\right)^{(N)}$ denote the corresponding vector space for this grading. If $A\in \left(H^0_{(\epsilon)}(\h^{\bot},d^{(\epsilon)}_{\h^{\bot},\mathfrak{q}})\right)^{(N)}$, let $A=\sum_{i=0}^N \epsilon^iF^{(N-i)}$. Then $d^{(\epsilon)}_{\h^{\bot},\mathfrak{q}}(A)=0\Leftrightarrow \epsilon d^{(1)}_{\h^{\bot},\mathfrak{q}}(F^{(N)})+\epsilon^2\left(d^{(2)}_{\h^{\bot},\mathfrak{q}}(F^{(N)})+d^{(1)}_{\h^{\bot},\mathfrak{q}}(F^{N-1)})\right)+\cdots =0\Leftrightarrow F=\sum_{i=0}^N F^{(N-i)}\in H^0(\h^{\bot},d_{\h^{\bot},\mathfrak{q}})$ by Prop. \ref{homogeneous}.

\noindent\textbf{ The specialization algebra} $\mathbf{H^0_{(\epsilon=1)}(\h^{\bot},d^{(\epsilon=1)}_{\h^{\bot},\mathfrak{q}}).}$ Denote as $<\epsilon-1>$ the ideal $(\epsilon-1)H^0_{(\epsilon)}(\h^{\bot},d^{(\epsilon)}_{\h^{\bot},\mathfrak{q}})$  of $H^0_{(\epsilon)}(\h^{\bot},d^{(\epsilon)}_{\h^{\bot},\mathfrak{q}})$ for the $\ast_{CF,\epsilon}-$ product. Define the specialized algebra  of $\h^{\bot}$ to be $H^0_{(\epsilon=1)}(\h^{\bot},d^{(\epsilon=1)}_{\h^{\bot},\mathfrak{q}}):=\left(H^0_{(\epsilon)}(\h^{\bot},d^{(\epsilon)}_{\h^{\bot},\mathfrak{q}})/<\epsilon-1>\right)$. Denote by $\ast_{CF,(\epsilon=1)}$ the Cattaneo-Felder product on $H^0_{(\epsilon=1)}(\h^{\bot},d^{(\epsilon=1)}_{\h^{\bot},\mathfrak{q}})$. Let $F=\sum_{i=0}^n F^{(n-i)},\;\;G=\sum_{j=0}^p G^{(p-j)}\in H^0(\h^{\bot},d_{\h^{\bot},\mathfrak{q}})$ be decompositions in homogeneous components (that is $\mathrm{deg}_{\mathfrak{q}}(F^{(k)})=k, \mathrm{deg}_{\mathfrak{q}}(G^{(m)})=m$) and let $F_{(\epsilon)}=\sum_{i=0}^n\epsilon^iF^{(n-i)},\;G_{(\epsilon)}=\sum_{j=0}^p\epsilon^jG^{(p-j)} \in H^0_{(\epsilon)}(\h^{\bot},d^{(\epsilon)}_{\h^{\bot},\mathfrak{q}})$. Consider the map $\mathfrak{i}_{\epsilon}: \; H^0(\h^{\bot},d_{\h^{\bot},\mathfrak{q}})\longrightarrow H^0_{(\epsilon)}(\h^{\bot},d^{(\epsilon)}_{\h^{\bot},\mathfrak{q}}),\;\;F\mapsto \mathfrak{i}_{\epsilon}(F)$. For $p<n$ we have $\mathfrak{i}_{\epsilon}(F+G)=F_{(\epsilon)}+ \epsilon ^{n-p}G_{(\epsilon)}=F_{(\epsilon)}+G_{(\epsilon)}+(\epsilon^{n-p}-~1)G_{(\epsilon)}$. Furthermore, the $\ast_{CF,\epsilon}$ product here is homogeneous of total degree $0$, that means for $A,B$ homogeneous in $\mathrm{deg}_{\mathfrak{q},\epsilon}$,  we have $\mathrm{deg}_{\mathfrak{q},\epsilon}\left(A\ast_{CF,\epsilon} B\right)=\mathrm{deg}_{\mathfrak{q},\epsilon}(A)+\mathrm{deg}_{\mathfrak{q},\epsilon}(B)$. Indeed, let $A=\sum_{j=0}^n \epsilon^{n-j}A^{(j)},\;B=\sum_{k=0}^p \epsilon^{p-k}B^{(k)}$ with $\mathrm{deg}_{\mathfrak{q},\epsilon}(A)=n, \mathrm{deg}_{\mathfrak{q},\epsilon}(B)=p$. Then $A\ast_{CF,\epsilon} B=AB+\sum_{i=1}^{\infty}\epsilon^i\sum_{\Gamma\in\mathbf{Q}_{i,2}^{(2)}}\overline{\omega}_{\Gamma}B_{\Gamma}(A,B)$. If for example $\Gamma\in \mathbf{Q}^{(2)}_{i,2}$ has one type I vertex that is not derived (i.e a root), then $\mathrm{deg}_{\mathfrak{q}}\left(B_{\Gamma}(A^{(j)},B^{(k)})\right)= j+k-(i+1)+1$. So $\forall i,\;\mathrm{deg}_{\mathfrak{q},\epsilon}\left(\epsilon^iB_{\Gamma}(\epsilon^{n-j}A^{(j)},\epsilon^{(p-k)}B^{(k)})\right)=n+p$. Thus $\mathfrak{i}_{\epsilon}(F\ast_{CF}G)=\mathfrak{i}_{\epsilon}(F)\ast_{CF,\epsilon}\mathfrak{i}_{\epsilon}(G)$. Let now $\pi_{(\epsilon=1)}:\;H^0_{(\epsilon)}(\h^{\bot},d^{(\epsilon)}_{\h^{\bot},\mathfrak{q}})\longrightarrow H^0_{(\epsilon=1)}(\h^{\bot},d^{(\epsilon=1)}_{\h^{\bot},\mathfrak{q}})$ be the canonical projection.
\newtheorem{s}[con]{Lemma}
\begin{s}\label{ftou}
Let $\g$ be a Lie algebra, $\h\subset\g$ a subalgebra and $\mathfrak{q}$ a subspace such that $\g=\h\oplus\mathfrak{q}$. The map $\mathfrak{i}_{(\epsilon=1)}:= \pi_{(\epsilon=1)}\circ\mathfrak{i}_{\epsilon}:\;H^0(\h^{\bot},d_{\h^{\bot},\mathfrak{q}})\longrightarrow H^0_{(\epsilon=1)}(\h^{\bot},d^{(\epsilon=1)}_{\h^{\bot},\mathfrak{q}})$ is an algebra isomorphism.
\end{s}

\noindent\textit{Proof.}  Let's first prove that $\mathfrak{i}_{(\epsilon=1)}$ is a morphism. For $F,G\in H^0(\h^{\bot},d_{\h^{\bot},\mathfrak{q}})$ and $p<n$ it is $ \mathfrak{i}_{(\epsilon=1)}(F+G)=\pi_{(\epsilon=1)}(F_{(\epsilon)}+\epsilon^{n-p}G_{(\epsilon)})=\pi_{(\epsilon=1)}(F_{(\epsilon)})+\pi_{(\epsilon=1)}(G_{(\epsilon)})+\pi_{(\epsilon=1)}((\epsilon^{n-p}-1)G_{(\epsilon)})=\mathfrak{i}_{(\epsilon=1)}(F)+\mathfrak{i}_{(\epsilon=1)}(G).$ Also by definition of the product on the quotient $H^0_{(\epsilon=1)}(\h^{\bot},d^{(\epsilon=1)}_{\h^{\bot},\mathfrak{q}})$, we have 
$\mathfrak{i}_{(\epsilon=1)}(F\ast_{CF}G)=\pi_{(\epsilon=1)}(F_{(\epsilon)}\ast_{CF,\epsilon} G_{(\epsilon)})=\pi_{(\epsilon=1)}(F_{(\epsilon)})\ast_{CF,(\epsilon=1)}\pi_{(\epsilon=1)}(G_{(\epsilon)}).$
In addition, $\mathfrak{i}_{(\epsilon=1)}$ is surjective because $\mathfrak{i}_{\epsilon}$ is surjective. Indeed, by the corollary after Prop. \ref{homogeneous}, $H^0_{(\epsilon)}(\h^{\bot},d^{(\epsilon)}_{\h^{\bot},\mathfrak{q}})$ is homogeneous with respect to the total degree $\mathrm{deg}_{\mathfrak{q},\epsilon}$. So an element $A$ of total degree $N$ can be written as $A=\sum_{i=0}^N\epsilon^i F^{(N-i)}$. Then $F:=\sum_{i=0}^NF^{(N-i)}\in H^0(\h^{\bot},d_{\h^{\bot},\mathfrak{q}})$ and $\mathfrak{i}_{\epsilon}(F)=A$. Finally $\mathfrak{i}_{(\epsilon=1)}$ is injective: For $F=\sum_{i=0}^NF^{(N-i)}\in H^0(\h^{\bot},d_{\h^{\bot},\mathfrak{q}})$ we have $\mathfrak{i}_{(\epsilon=1)}(F)=0\Leftrightarrow \mathfrak{i}_{\epsilon}(F)\in (\epsilon-1)H^0_{(\epsilon)}(\h^{\bot},d^{(\epsilon)}_{\h^{\bot},\mathfrak{q}})$. That is, $F=\sum_{i=0}^N\epsilon^i F^{(N-i)}\in (\epsilon-1)H^0_{(\epsilon)}(\h^{\bot},d^{(\epsilon)}_{\h^{\bot},\mathfrak{q}})$ by the discussion before this Lemma. Then for $\epsilon=1$ we take $ \sum_{i=0}^NF^{(N-i)}=0$ and thus $F=0$. $\diamond$

\noindent\textbf{3.3 The affine case: The algebras} $\mathbf{H^0_{(\epsilon)}(\h_{\lambda}^{\bot},d^{(\epsilon)}_{\h^{\bot}_{\lambda},\mathfrak{q}})}$ \textbf{and} $\mathbf{H^0_{(\epsilon=1)}(\h_{\lambda}^{\bot},d^{(\epsilon=1)}_{\h^{\bot}_{\lambda},\mathfrak{q}})}$. Set $\h_{\lambda}$ to be the vector subspace of $S(\g)$ generated by the set $\{H+\lambda(H),\;H\in\h\}$ . Let $\h_{\lambda}^{\bot}:=\{f\in\g^{\ast}/f|_{\h}=-\lambda\}$ an affine subspace of $\g^{\ast}$. We will  abusively write $-\lambda+\h^{\bot}=\h_{\lambda}^{\bot}$ since $\lambda$ is defined only on $\h$. Define the reduction algebra $H^0_{(\epsilon)}(\h_{\lambda}^{\bot},d^{(\epsilon)}_{\h^{\bot}_{\lambda},\mathfrak{q}})$ using homogeinity with respect to the deformation parameter $\epsilon$ in the same way as in the vector space case of $\mathcal{x}$ 3.2. This is less easy now: the root of a $\mathcal{B}-$ type graph can be either in $\mathfrak{q}^{\ast}$ or $\h^{\ast}$ as before. If the root is an $H\in\h^{\ast}$ then the existence of the character $\lambda$ would give $-\lambda(H)\in\mathbb{R}$ after evaluation at the root. So this time the operators $B_{\Gamma}$  are not homogeneous with respect to $\mathrm{deg}_{\mathfrak{q}}$ and consequently the  equations $d^{(\epsilon)}_{\h^{\bot}_{\lambda},\mathfrak{q}}(F_{(\epsilon)})=0$ are not homogeneous with respect to $\mathrm{deg}_{\mathfrak{q},\epsilon}$.\\
\noindent Denote by $H^0_{(\epsilon)}(\h_{\lambda}^{\bot},d^{(\epsilon)}_{\h^{\bot}_{\lambda},\mathfrak{q}})$ the algebra of polynomials $P_{(\epsilon)}$, solutions of the equation $d^{(\epsilon)}_{\h^{\bot}_{\lambda},\mathfrak{q}}(P_{(\epsilon)})=0$, equipped with the $\ast_{CF,\epsilon}$ product. It will be called the \textsl{reduction algebra} over $-\lambda+\h^\bot$. Denote by $H^0_{(\epsilon=1)}(\h_{\lambda}^{\bot},d^{(\epsilon=1)}_{\h_{\lambda}^{\bot},\mathfrak{q}})$ the \textsl{specialized reduction algebra} $H^0_{(\epsilon)}(\h_{\lambda}^{\bot},d^{(\epsilon)}_{\h_{\lambda}^{\bot},\mathfrak{q}})/<\epsilon-1>$ over $-\lambda+\h^\bot$ with the corresponding  Cattaneo-Felder product denoted by $\ast_{CF,(\epsilon=1)}$. Let $d_{\h_{\lambda}^{\bot},\mathfrak{q}}:\;S(\mathfrak{q})\longrightarrow S(\mathfrak{q})\otimes \h$ be the differential operator $d_{\h_{\lambda}^{\bot},\mathfrak{q}}= \sum_{i=1}^{\infty}d^{(i)}_{\h_{\lambda}^{\bot},\mathfrak{q}}$ where $d^{(i)}_{\h_{\lambda}^{\bot},\mathfrak{q}}=\sum_{\substack{\Gamma\in \mathcal{B}_i\cup \mathcal{BW}_i}}\overline{\omega}_{\Gamma}B_{\Gamma}$, and let $H^0(\h_{\lambda}^{\bot},d_{\h^{\bot}_{\lambda},\mathfrak{q}})$ be the vector space of polynomials $P$,  solutions of the equation $d_{\h^{\bot}_{\lambda},\mathfrak{q}}(P)=0$. The Cattaneo-Felder construction without $\epsilon$ is still valid for functions on $-\lambda+\h^{\bot}$ and defines an associative product $\ast_{CF}$ on $H^0(\h_{\lambda}^{\bot},d_{\h_{\lambda}^{\bot},\mathfrak{q}})$. In the sequence $H^0(\h_{\lambda}^{\bot},d_{\h_{\lambda}^{\bot},\mathfrak{q}})$ will stand for the algebra $\left(H^0(\h_{\lambda}^{\bot},d_{\h_{\lambda}^{\bot},\mathfrak{q}}),\ast_{CF}\right)$. In the affine case, $H^0_{(\epsilon=1)}(\h_{\lambda}^{\bot},d^{(\epsilon=1)}_{\h^{\bot}_{\lambda},\mathfrak{q}})$ and $H^0(\h_{\lambda}^{\bot},d_{\h^{\bot}_{\lambda},\mathfrak{q}})$ are not isomorphic as in the vector space case. However the following is true: Let $F^{'}=\sum_{i=0}^q\epsilon^iF^{'}_i\in H^0_{(\epsilon)}(\h_{\lambda}^{\bot},d^{(\epsilon)}_{\h^{\bot}_{\lambda},\mathfrak{q}})$, and consider the linear map $J:\;H^0_{(\epsilon)}(\h_{\lambda}^{\bot},d^{(\epsilon)}_{\h^{\bot}_{\lambda},\mathfrak{q}})\longrightarrow H^0(\h_{\lambda}^{\bot},d_{\h^{\bot}_{\lambda},\mathfrak{q}})\;,\;\;J(F^{'})=\sum_k F^{'}_k.$ Obviously $J(<\epsilon-1>)=0$ and we denote as $\overline{J}:\;H^0_{(\epsilon=1)}(\h_{\lambda}^{\bot},d^{(\epsilon=1)}_{\h^{\bot}_{\lambda},\mathfrak{q}})\longrightarrow H^0(\h_{\lambda}^{\bot},d_{\h^
{\bot}_{\lambda},\mathfrak{q}})$ the induced quotient map. We borrow notation from $\mathcal{x}$ 3.2 and let $\pi_{(\epsilon=1)}:\;H^0_{(\epsilon)}(\h_{\lambda}^{\bot},d^{(\epsilon)}_{\h_{\lambda}^{\bot},\mathfrak{q}})\longrightarrow H^0_{(\epsilon=1)}(\h_{\lambda}^{\bot},d^{(\epsilon=1)}_{\h^{\bot}_{\lambda},\mathfrak{q}})$ denote here also the canonical projection.
\newtheorem{stiros}[con]{Lemma}
\begin{stiros}\label{jbar}
With these notations, $\overline{J}:\;H^0_{(\epsilon=1)}(\h_{\lambda}^{\bot},d^{(\epsilon=1)}_{\h^{\bot}_{\lambda},\mathfrak{q}})\longrightarrow H^0(\h_{\lambda}^{\bot},d_{\h^{\bot}_{\lambda},\mathfrak{q}})$ is injective.
\end{stiros}
\textit{Proof.}  The map $J$ is well defined. Indeed, by definition of $H^0_{(\epsilon)}(\h_{\lambda}^{\bot},d^{(\epsilon)}_{\h^{\bot}_{\lambda},\mathfrak{q}})$, if $F^{'}=\sum_{i=0}^q\epsilon^iF^{'}_i\in H^0_{(\epsilon)}(\h_{\lambda}^{\bot},d^{(\epsilon)}_{\h^{\bot}_{\lambda},\mathfrak{q}})$ (the terms $F^{'}_i$ are not homogeneous polynomials here) then it satisfies the equations $d_{\h_{\lambda}^{\bot},\mathfrak{q}}^{(1)}(F^{'}_0)=0,\;\;d_{\h_{\lambda}^{\bot},\mathfrak{q}}^{(1)}(F^{'}_{1})+d_{\h_{\lambda}^{\bot},\mathfrak{q}}^{(2)}(F^{'}_0)=0$ and etc. Thus $\left(\sum_{p\geq 1}\epsilon^p d^{(p)}_{\h_{\lambda}^{\bot},\mathfrak{q}}\right)\left(\sum_{k\geq 0}\epsilon^kF_k^{'}\right)=0$ and so 
\[J\left(\left(\sum_{p\geq 1}\epsilon^p d^{(p)}_{\h_{\lambda}^{\bot},\mathfrak{q}}\right)\left(\sum_{k\geq 0}\epsilon^kF_k^{'}\right)\right)=\left(\sum_{p\geq 1}d^{(p)}_{\h_{\lambda}^{\bot},\mathfrak{q}}\right)\left(\sum_{k\geq 0}F_k^{'}\right)=0\]
 Thus $J(F^{'})=\sum_{k\geq 0}F_k^{'}\in  H^0(\h_{\lambda}^{\bot},d_{\h^{\bot}_{\lambda},\mathfrak{q}})$. Then $J(F^{'})=0\Rightarrow F^{'}=\sum_{k\geq 1}\epsilon^kF^{'}_k-J(F^{'})=\sum_{k\geq 1}\epsilon^kF^{'}_k-\sum_{k\geq 1}F^{'}_k=\sum_{k\geq 1}(\epsilon^k-1)F^{'}_k\in<\epsilon-1>,$
that is  $J(F^{'})=0\Rightarrow F'\in<\epsilon-1>$. So $F'=0$ in the quotient $H^0_{(\epsilon=1)}(\h_{\lambda}^{\bot},d^{(\epsilon=1)}_{\h^{\bot}_{\lambda},\mathfrak{q}})$. Thus $\overline{J}$ is injective. It's clearly a map of algebras for the products $\ast_{CF,(\epsilon=1)}$ and $\ast_{CF}$.

\section{The biquantization diagram of  $\g^{\ast},-\lambda+\h^{\bot}$.}
\textbf{4.1 Biquantization diagrams and the bimodule structure.} We briefly return to the general setting of a Poisson manifold $X$.  The construction of $\mathcal{x}$ 2.3 for essentially  one coisotropic submanifold since we considered $C_1=X,\;C_2\subset X$ and thus only two colors, is generalized for two coisotropic submanifolds $C_1,C_2\subset X$ and four colors as follows. Consider the diagram of Fig. 2 where on each semi-axis we have the one coisotropic submanifold case of $\mathcal{x}$ 2.3 with respect to the computation of the coefficients $\overline{\omega}_{\Gamma}$ and the differential operators $B_{\Gamma}$.  Let $\mathcal{H}^{++}=\{z\in\mathbb{C}/ \Im(z)>0,\Re(z)>0\}$ and $\epsilon_1,\epsilon_2$ take values in $\{-,+\}$. For $z_1,z_2\in \mathcal{H}^{++}$ define the angle function of four colors as $\phi_{\epsilon_1,\epsilon_2}(z_1,z_2):=\mathrm{arg}(z_1-z_2) +\epsilon_1\mathrm{arg}(z_1-\overline{z_2})+\epsilon_2\mathrm{arg}(z_1+~\overline{z_2})+\epsilon_1\epsilon_2\mathrm{arg}(z_1+z_2)$. The angle form $\mathrm{d}\phi_{\epsilon_1\epsilon_2}$ is defined correspondingly. When we consider concentrations of points close to the vertical or the horizontal axis, the angle form of four colors $\mathrm{d}\phi_{\epsilon_1,\epsilon_2}$ degenerates to the angle form of two colors and the situation trivializes to the case of one coisotropic submanifold. For example in Fig. 2, if we concentrate at the vertical axis, an edge has two options of colors namely $+$ and $-$ (corresponding in the 4-color situation to $(+,-)$ and $(-,-)$) whether the variable assigned to that edge belongs or doesn't belong respectively to $C_1$. An important feature in ~\cite{CF2} is the existence of a bimodule structure for the reduction space at the origin of the biquantization diagram. Let $H^0_{(\epsilon)}(C_1,d^{(\epsilon)}_{C_1}),\;H^0_{(\epsilon)}(C_2,d^{(\epsilon)}_{C_2}),\;H^0_{(\epsilon)}(C_1\cap C_2,d^{(\epsilon)}_{C_1,C_2})$ be the reduction spaces at the vertical axis,  the horizontal axis and the corner of the diagram (corresponding to $C_1\cap C_2$) respectively. The first two are in fact algebras with the corresponding $\ast_{CF,\epsilon}$ product, while the third is not. In generality these spaces correspond to the 0-th cohomology of a Chevalley-Eilenberg complex, which is actually an $A_{\infty}-$ algebra in the case of the first two. The differential $d^{(\epsilon)}_{C_1,C_2}$ is defined in terms of graphs in two colors as in $\mathcal{x}$ 2.3 while for $d^{(\epsilon)}_{C_i}$  one uses a 4-colored labelling function $L$.
\begin{figure}[h!]
\begin{center}
\includegraphics[width=8cm]{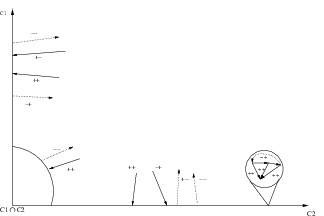}
\caption{\footnotesize The biquantization diagram for $C_1,\;C_2$ with colors and concentrations.}
\end{center}
\end{figure}
The authors in ~\cite{CF2} construct two module structures $\ast_1:\; H^0_{(\epsilon)}(C_1,d^{(\epsilon)}_{C_1})\times H^0_{(\epsilon)}(C_1\cap C_2,d^{(\epsilon)}_{C_1,C_2})\longrightarrow H^0_{(\epsilon)}(C_1\cap C_2,d^{(\epsilon)}_{C_1,C_2}),\;\;<f,\rho>\mapsto f\ast_1\rho$, and $\ast_2: \; H^0_{(\epsilon)}(C_1\cap C_2,d^{(\epsilon)}_{C_1,C_2})\times H^0_{(\epsilon)}(C_2,d^{(\epsilon)}_{C_2})\longrightarrow H^0_{(\epsilon)}(C_1\cap C_2,d^{(\epsilon)}_{C_1,C_2}),\;
<\rho,g>\mapsto \rho\ast_2g$.
By the original construction, we know exactly the formulas for these module structures. If we denote by $\mathbf{Q}_{k,2}^{(4)}$ the family of admissible graphs of four colors and two type II vertices, then for $\rho \in H^0_{(\epsilon)}(C_1\cap~ C_2,d^{(\epsilon)}_{C_1,C_2}),\; f\in H^0_{(\epsilon)}(C_1,d^{(\epsilon)}_{C_1})$ one writes $f \ast_1 \rho =f\cdot\rho + \sum_{k=1}^{\infty} \frac{\epsilon ^{k}}{k!}\sum_{\Gamma \in \mathbf{Q}_{k,2}^{(4)}}\overline{\omega}_{\Gamma}B_{\Gamma}(f,\rho)$,
for the left $ H^0_{(\epsilon)}(C_1,d^{(\epsilon)}_{C_1})-$ module structure of $H^0_{(\epsilon)}(C_1\cap~ C_2,d^{(\epsilon)}_{C_1,C_2})$ and similarly, for the horizontal axis, the right $H^0_{(\epsilon)}(C_2,d^{(\epsilon)}_{C_2})-$ module structure of $H^0_{(\epsilon)}(C_1\cap~ C_2,d^{(\epsilon)}_{C_1,C_2})$ is written as $\rho \ast_2 g =\rho\cdot g + \sum_{k=1}^{\infty} \frac{\epsilon ^{k}}{k!}\sum_{\Gamma\in \mathbf{Q}_{k,2}^{(4)}}\overline{\omega}_{\Gamma}B_{\Gamma}(\rho,g)$. Note that here $B_{\Gamma}$ is defined as in $\mathcal{x}$ 2.2  but with the use of a 4-colored labelling function $L$.

\noindent\textbf{4.2 The operators} $\mathbf{T_1}$ \textbf{and} $\mathbf{T_2.}$ Let $q(Y) := \det_{\g} \left(\frac{\sinh\frac{\mathrm{ad} Y}{2}}{\frac{adY}{2}}\right), Y\in\g$ and $\beta:\;S(\g)\longrightarrow U(\g)$ be the PBW symmetrizaton map. We use the notation $\ast_D$ for the star-product on $S(\g)$ written without the parameter $\epsilon$ and $\ast_{DK}$ for the star-product written with $\epsilon$. Recall from \cite{K} $\mathcal{x}$ 8.3, Theorem 8.2, that the Kontsevich product $\ast_D$, satisfies the relation $ \beta\left(\partial_{q^{\frac{1}{2}}}(f\ast_D g)\right)=\beta(\partial_{q^{\frac{1}{2}}}f)\cdot \beta(\partial_{q^{\frac{1}{2}}}g)$.  Consider now the biquantization diagram as in Fig. 2, putting the coisotropic submanifold $\mathfrak{o}^\bot=\g^{\ast}$ on the vertical axis  and $-\lambda+\h^{\bot}$ on the horizontal one. Let $\mathfrak{q}$ be such that $\g=\h\oplus\mathfrak{q}$, $\Gamma\in\mathbf{Q}^{(4)}_{n,2}$ and $c_e=(\cdot,\cdot)$ denote the color of an edge $e\in E(\Gamma)$.  Let again $\{H_1^{\ast},\ldots,H_t^{\ast},Q_1^{\ast},\ldots,Q_r^{\ast}\}$ be a basis of $\g^{\ast}$ formed from bases of $\h$ and $\mathfrak{q}$ respectively and let $(x_i)_{i=1,\ldots,n}$ be a coordinate system relatively to this basis.  We assign a variable $x_{s,e}\in\{x_1,\ldots,x_t\}$ if $e\in E(\Gamma)$ is of color $c_e=(\pm,-)$ and a variable $x_{k,e}\in\{x_{t+1},\ldots,x_{t+r}\}$ if $c_e=\;(\pm,+).$
Since there is no variable in $\g^{\bot}$, there are only two colors in this diagram namely $(+,-)$ (corresponding to $\h^\ast),$ or $(+,+)$ (corresponding to $\mathfrak{q}^\ast\simeq\h^{\bot}).$ Graphs here may have double edges and still satisfy $\overline{\omega}_{\Gamma}\neq 0$ as long as these edges are of different color. \\
\noindent Let $H^0_{(\epsilon)}(\h_{\lambda}^{\bot},d^{(\epsilon)}_{\g^{\ast},\h_{\lambda}^{\bot},\mathfrak{q}})$ be the reduction space at the corner of this diagram. We use the bimodule structure of $\mathcal{x}$ 4.1: Let $\ast_1$ be the left $H^0_{(\epsilon)}(\g^{\ast},d^{(\epsilon)}_{\g^{\ast}})-$ module structure of $H^0_{(\epsilon)}(\h_{\lambda}^{\bot},d^{(\epsilon)}_{\g^{\ast},\h_{\lambda}^{\bot},\mathfrak{q}})$ and consider the map $T_1:H^0_{(\epsilon)}(\g^{\ast},d^{(\epsilon)}_{\g^{\ast}})\longrightarrow H^0_{(\epsilon)}(\h_{\lambda}^{\bot},d^{(\epsilon)}_{\g^{\ast},\h_{\lambda}^{\bot},\mathfrak{q}}),\; F \mapsto F\ast_1 1$. Similarly for the horizontal axis let $\ast_2$ be the right $H^0_{(\epsilon)}(\h^{\bot}_{\lambda},d^{(\epsilon)}_{\h^{\bot}_{\lambda},\mathfrak{q}}) -$ module structure of $H^0_{(\epsilon)}(\h_{\lambda}^{\bot},d^{(\epsilon)}_{\g^{\ast},\h_{\lambda}^{\bot},\mathfrak{q}})$ and let $T_2:H^0_{(\epsilon)}(\h^{\bot}_{\lambda},d^{(\epsilon)}_{\h^{\bot}_{\lambda},\mathfrak{q}}) \longrightarrow H^0_{(\epsilon)}(\h_{\lambda}^{\bot},d^{(\epsilon)}_{\g^{\ast},\h_{\lambda}^{\bot},\mathfrak{q}})$, $G\mapsto 1 \ast_2 G$.
The calculation of the reduction algebras $H^0_{(\epsilon)}(\h_{\lambda}^{\bot},d^{(\epsilon)}_{\h_{\lambda}^{\bot},\mathfrak{q}})$, $H^0_{(\epsilon)}(\g^{\ast},d^{(\epsilon)}_{\g^{\ast}})$ and the reduction space $H^0_{(\epsilon)}(\h_{\lambda}^{\bot},d^{(\epsilon)}_{\g^{\ast},\h_{\lambda}^{\bot},\mathfrak{q}})$ is done in ~\cite{CT}. Recall that $H^0_{(\epsilon)}(\g^{\ast},d^{(\epsilon)}_{\g^{\ast}})\simeq U_{(\epsilon)}(\g)$ as associative algebras and $H^0_{(\epsilon)}(\h_{\lambda}^{\bot},d^{(\epsilon)}_{\g^{\ast},\h_{\lambda}^{\bot},\mathfrak{q}})\simeq S(\mathfrak{q})[\epsilon]$ as vector spaces, so we will write $T_1:\left(S_{(\epsilon)}(\g),\ast_{DK}\right)\simeq \left(U_{(\epsilon)}(\g),\cdot\right)\longrightarrow S(\mathfrak{q})[\epsilon]$,  and $T_2:\;H^0_{(\epsilon)}(\h^{\bot}_{\lambda},d^{(\epsilon)}_{\h^{\bot}_{\lambda},\mathfrak{q}})\longrightarrow S(\mathfrak{q})[\epsilon]$. Let $F\in (U_{(\epsilon)}(\g),\cdot),\;\; G \in H^0_{(\epsilon)}(\h^{\bot}_{\lambda},d^{(\epsilon)}_{\h^{\bot}_{\lambda},\mathfrak{q}})$. From \cite{CF2}, the compatibility relation between the two module structures is 
\begin{equation}\label{compatibility}
(F\ast_1 1)\ast_2 G= G\ast_1 (1\ast_2 G),
\end{equation}
For $F\in S(\g)$, a possible graph in $T_1(F)$ is not necessarily of $\mathcal{W}-$ type because $F$ can be derived by two different colors, $(+,-)$ and $(+,+)$. For example the graph $\Gamma\in\mathbf{Q}^{
(4)}_{1,1}$ with two edges towards $F$ colored by $(+,-)$ and $(+,+)$ respectively, contributes non-trivially to $F\ast_1 1$. However if $F\in S(\mathfrak{q})$, then $F$ can receive only $(+,+)\;$- colored edges. Thus the previous graph has now a double edge of the same color and it is not allowed. In this case, $T_1$ is composed only of $\mathcal{W}-$ type graphs, with color $(+,+)$ for the edges deriving $F$, and arbitrary colors for the edges in the wheel.

\noindent \textbf{4.3 Decompositions and projections.} Consider now the graph $\Gamma'$ of Fig. 3 in the biquantization diagram of $\g^\ast$, $-\lambda+\h^\bot$. According to \cite{CRT} this graph is a \textsl{small loop} and it is allowed for $T_1$. Furthermore, the authors in \cite{CRT} explain how to compute its coefficient $\overline{\omega}_{\Gamma'}$. One has to use the angle form $\mathrm{d}arg(z)$ with $z$ being the position of the type I vertex of the graph, that is $\overline{\omega}_{\Gamma'}=\frac{1}{2}$. Furthermore, the corresponding operator $B_{\Gamma'}$ is the trace of the $\ad-$ action. For $H\in \h$ of degree one as a polynomial function, it is easy to see with respect to $\mathcal{x}$ 2.1 that $\rho(H):=-\overline{\omega}_{\Gamma'}B_{\Gamma'}(H)$. We need two decompositions which we recall from \cite{CT}. Let $\overline{\beta}_{\mathfrak{q}}:\;S(\mathfrak{q})\longrightarrow~(U(\g)/U(\g)\h_{\lambda+\rho})$ be the quotient symmetrization map. Then one has the decomposition $
U(\g)=\overline{\beta}_{\mathfrak{q}}\circ\partial_{q^{\frac{1}{2}}}(S(\mathfrak{q}))\oplus U(\g)\cdot\h_{\lambda+\rho}$ \footnote{The second summand is just $U(\g)\h_{\lambda+\rho}$. We make the product visible for this ideal because right afterwards we consider $\ast-$ ideals.}. Applying the isomorphism $(S(\g),\ast_{D})\stackrel{\simeq}{\longrightarrow}(U(\g),\cdot)$, we get the decomposition $S(\g)=S(\mathfrak{q})\oplus S(\g)\ast_{D}\h_{\lambda+\rho}$. Analogous decompositions hold for $S_{(\epsilon)}(\g),U_{(\epsilon)}(\g)$ as defined at the end of $\mathcal{x}$ 2.1. Indeed,
\begin{equation}\label{decomposition S(g)[[e]]}
S_{(\epsilon)}(\g)=S(\mathfrak{q})[\epsilon]\oplus S_{(\epsilon)}(\g)\ast_{DK}\h_{\lambda+\rho}
\end{equation}
The PBW theorem holds also for the algebras $S_{(\epsilon)}(\g),\; U_{(\epsilon)}(\g)$ and so there is a symmetrization map $\beta_{(\epsilon)}:\;S_{(\epsilon)}(\g)\longrightarrow~U_{(\epsilon)}(\g)$. We denote as $\overline{\beta}_{\mathfrak{q},(\epsilon)}:\;S(\mathfrak{q})[\epsilon]\longrightarrow U_{(\epsilon)}(\g)/U_{(\epsilon)}(\g)\h_{\lambda+\rho}$ the quotient of this symmetrization map with respect to $\mathfrak{q}$. Set now $q_{(\epsilon)}(X):=q(\epsilon X)$ for $X\in\g$. Using the isomorphism $(S_{(\epsilon)}(\g),\ast_{DK})\simeq (S_{(\epsilon)}(\g),\ast_{CF})\simeq(U_{(\epsilon)}(\g),\cdot)$, the deformed algebra
 $U_{(\epsilon)}(\g)$ can be decomposed as
\begin{equation}\label{decomposition U(g)[[e]]}
U_{(\epsilon)}(\g)=\overline{\beta}_{\mathfrak{q},(\epsilon)}\circ\partial_{q_{(\epsilon)}^{\frac{1}{2}}}(S(\mathfrak{q})[\epsilon])\oplus U_{(\epsilon)}(\g)\cdot\h_{\lambda+\rho}.
\end{equation}

\newtheorem{sti}[con]{Lemma}
\begin{sti}\label{lemma for H}
Consider the biquantization diagram with $\g^{\ast}$ on the vertical, and $-\lambda+\h^{\bot}$ on the horizontal axis.  Let $H\in \h$ be a function of degree 1 on the vertical axis, and let $\ast_1$ denote the left $(S_{(\epsilon)}(\g),\ast_{DK})-$ module structure of  $H^0_{(\epsilon)}(\h_{\lambda}^{\bot},d^{(\epsilon)}_{\g^{\ast},\h_{\lambda}^{\bot},\mathfrak{q}})$. Then
\begin{equation}\label{Hast1=0}
(H+(\lambda+\rho)(H))\ast_1 1=0.
\end{equation}
\end{sti}
\textit{Proof.} As $H$ is a function of degree one, it can receive only one edge colored by $(+,-)$ and so does $H+(\lambda+\rho)(H)$. At the corner of the diagram, derivation of $1$ will give a zero outcome so we look for graphs deriving only the function $H+(\lambda+\rho)(H)$ . Let $k$ be the number of type I vertices in a possible graph $\Gamma$ in the expression $(H+(\lambda+\rho)(H))\ast_1 1$. Then from $\mathcal{x}$ 2.2, $\Gamma$ should have $2k$ edges and satisfy the restriction $2k-1 \leq k$. So $k\leq1$ and the only possible graph here is $\Gamma'$. By the end of $\mathcal{x}$ 4.1, $(H+(\lambda+\rho)(H))\ast_1 1=(H+(\lambda+\rho)(H))\cdot~1+\overline{\omega}_{\Gamma'} B_{\Gamma'}((H+(\lambda+\rho)(H)),1)$. For the second summand we have $B_{\Gamma'}(H+(\lambda+\rho)(H),1)=\mathrm{Tr}\;\mathrm{ad}_\g(H)$ by definition of this differential operator. Thus $(H+(\lambda+\rho)(H))\ast_1 1=(H+(\lambda+\rho)(H))\cdot 1 +\overline{\omega}_{\Gamma'}\mathrm{Tr}\;\mathrm{ad}_\g(H)$. Since we finally have to restrict this result to $-\lambda+\h^{\bot}$, one gets $-\lambda[(H+\lambda(H))\ast_1 1]=-\lambda[(H+(\lambda+\rho)(H))\cdot 1] +\overline{\omega}_{\Gamma'}\mathrm{Tr}\; \mathrm{ad}_\g(H)=0$. $\diamond$  
\begin{figure}[h!]
\begin{center}
\includegraphics[width=5cm]{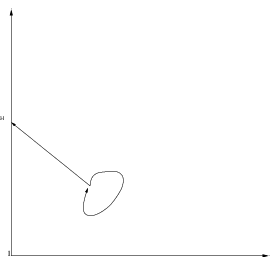}
\caption{\footnotesize The graph $\Gamma'$. The only possible graph in $(H+(\lambda+\rho)(H))\ast_1 1$.}
\end{center}
\end{figure}

The operator $T_1(F)=F\ast_1 1$ is not a priori of constant coefficients. This is why we are now allowed to have double edges deriving $F$ in the graphs contributing to $T_1$. For such a graph $\Gamma$, we have $\overline{\omega}_{\Gamma}\neq 0$ in general as the double edges might be of different color e.g (+,+)(for an edge carrying a variable in $\mathfrak{q}^{\ast}$) and (+,-) (for an edge carrying a variable in $\h^{\ast}$). This was not allowed on the horizontal axis for the algebra $H^0_{(\epsilon)}(\h^{\bot}_{\lambda},d^{(\epsilon)}_{\h^{\bot}_{\lambda},\mathfrak{q}})$. The following uses arguments of \cite{CT} and should give better insight on $T_1$.

\newtheorem{spahia}[con]{Lemma}
\begin{spahia}\label{lemma for T_1}
The operator $T_1$ in the biquantization diagram of $\g^{\ast}$ and $-\lambda+\h^{\bot}$ satisfies
\begin{equation}\label{T_1 is projection}
F\ast_1 1=T_1(F)=0\Longleftrightarrow F\in S_{(\epsilon)}(\g)\ast_{DK}\h_{\lambda+\rho}.
\end{equation}
\end{spahia}
\textit{Proof.} Let $F\in H^0_{(\epsilon)}(\g^{\ast},d^{(\epsilon)}_{\g^{\ast}})\simeq U_{(\epsilon)}(\g)\simeq (S_{(\epsilon)}(\g),\ast_{DK})$ be at the vertical axis. 
By the decomposition (\ref{decomposition S(g)[[e]]}), let $F=A+B$ with $A\in S(\mathfrak{q})[\epsilon],\;B\in S_{(\epsilon)}(\g)\ast_{DK}\h_{\lambda+\rho}$. Then $T_1(F)=F\ast_1 1=A\ast_1 1=T_1(A)$. Since $A\in S(\mathfrak{q})[\epsilon]$, the operator $T_1$ acting on $A$ is composed only of $\mathcal{W}-$ type graphs of one color. That is $T_1(A)=\exp\left(\sum_i(\sum_{W^i\in\mathcal{W}_i}\overline{\omega}_{W^i}B_{W^i})\right)(A)$. So in this case, $T_1$ is of constant coefficients and invertible because it is of exponential type. This means that $T_1(F)=0\Rightarrow T_1(A)=0\Rightarrow A=0$ and $F=B\in S_{(\epsilon)}(\g)\ast_D\h_{\lambda+\rho}$. Conversely, let $Q\ast_{DK} (H+(\lambda+\rho)(H))\in S_{(\epsilon)}(\g)\ast_{DK}\h_{\lambda+\rho}$. Then $T_1\left(Q\ast_{DK}(H+(\lambda+\rho)(H))\right)=\left(Q\ast_{DK}(H+(\lambda+\rho)(H))\right)\ast_1 1=Q\ast_1\left((H+(\lambda+\rho)(H))\ast_1 1\right)=0$ where for the last equation we used Lemma \ref{lemma for H} and the second, the compatibility property (\ref{compatibility}) of the $\ast$ product. $\diamond$

Let $\overline{T}_1:=T_1|_{S(\mathfrak{q})[\epsilon]}$. By lemma \ref{lemma for T_1}, $\overline{T}_1$ is a vector space isomorphism and we will write abusively its inverse $\overline{T}_1^{-1}:\;S(\mathfrak{q})[\epsilon]\longrightarrow S(\mathfrak{q})[\epsilon]\subset H^0_{(\epsilon)}(\g^{\ast},d^{(\epsilon)}_{\g^{\ast}})$. The operators $\overline{T}_1, T_2$ are constant coefficient operators described by $\mathcal{W}-$ type graphs. From now on we identify the products $\ast_{DK}$ and $\ast_{CF,\epsilon}$ on $S_{(\epsilon)}(\g)$.

\section{$H^0_{(\epsilon)}(\h^{\bot}_{\lambda},d^{(\epsilon)}_{\h^{\bot}_{\lambda},\mathfrak{q}})\simeq (U_{(\epsilon)}(\g)/U_{(\epsilon)}(\g)\h_{\lambda})^{\h}$.}
Recall from \cite{CF2}, \cite{CT} that the map $\overline{\beta}_{\mathfrak{q},(\epsilon)}\circ\partial_{q_{(\epsilon)}^{\frac{1}{2}}} \circ \overline{T}_1^{-1}T_2:\;\; (H^0_{(\epsilon)}(\h^{\bot}_{\lambda},d^{(\epsilon)}_{\h^{\bot}_{\lambda},\mathfrak{q}}),\ast_{CF})\longrightarrow\left(U_{(\epsilon)}(\g)/U_{(\epsilon)}(\g)\h_{\lambda}\right)$,
satisfies $mod[U_{(\epsilon)}(\g)\h_{\lambda}]$ the relation
\[\left(\overline{\beta}_{\mathfrak{q},(\epsilon)}\circ\partial_{q_{(\epsilon)}^{\frac{1}{2}}}\circ \overline{T}_1^{-1}T_2\right)(F_1\ast_{CF}F_2)\equiv \left[\left(\overline{\beta}_{\mathfrak{q},(\epsilon)}\circ\partial_{q_{(\epsilon)}^{\frac{1}{2}}}\circ \overline{T}_1^{-1}T_2\right)(F_1)\cdot\left(\overline{\beta}_{\mathfrak{q},(\epsilon)}\circ \partial_{q_{(\epsilon)}^{\frac{1}{2}}}\circ \overline{T}_1^{-1}T_2\right)(F_2)\right].\] 
\newtheorem{oki}[con]{Theorem}
\begin{oki}\label{main}
Let $\g$ be a Lie algebra and $\h\subset \g$ a Lie subalgebra. Let also $\lambda\in \h^{\ast}$, $\lambda([\h,\h])=0$ be a character of $\h$ and choose a supplementary space $\mathfrak{q}$ for $\h$ in $\g$.  The map $\overline{\beta}_{\mathfrak{q},(\epsilon)}\circ\partial_{q_{(\epsilon)}^{\frac{1}{2}}}\circ \overline{T}_1^{-1}T_2$ is an non-canonical algebra isomorphism between $H^0_{(\epsilon)}(\h^{\bot}_{\lambda},d^{(\epsilon)}_{\h^{\bot}_{\lambda},\mathfrak{q}})$  and $\left(U_{(\epsilon)}(\g)/U_{(\epsilon)}(\g)\h_{\lambda+\rho}\right)^{\h}$,
\begin{equation}\label{main isomorphism}
\overline{\beta}_{\mathfrak{q},(\epsilon)}\circ\partial_{q_{(\epsilon)}^{\frac{1}{2}}}\circ \overline{T}_1^{-1}T_2:\; H^0_{(\epsilon)}(\h^{\bot}_{\lambda},d^{(\epsilon)}_{\h^{\bot}_{\lambda},\mathfrak{q}})\stackrel{\sim}{\longrightarrow} \left(U_{(\epsilon)}(\g)/U_{(\epsilon)}(\g)\h_{\lambda+\rho}\right)^{\h}.
 \end{equation}
\end{oki}
We give the proof of Theorem \ref{main} in two parts, through the next Propositions of this section.
\newtheorem{okae}[con]{Proposition}
\begin{okae}\label{direct main}
With these notations, $H^0_{(\epsilon)}(\h^{\bot}_{\lambda},d^{(\epsilon)}_{\h^{\bot}_{\lambda},\mathfrak{q}})\subset (U_{(\epsilon)}(\g)/U_{(\epsilon)}(\g)\h_{\lambda+\rho})^{\h}$.
\end{okae}
\textit{Proof}.  Let $f\in\g^{\ast}$ such that $f|_{\h}=\lambda$. Fix a biquantization diagram with $\g^{\ast}$ at the vertical axis, and $-f+\h^{\bot}$ at the horizontal. Let $F\in H^0_{(\epsilon)}(\h^{\bot}_{\lambda},d^{(\epsilon)}_{\h^{\bot}_{\lambda},\mathfrak{q}})$ be on the horizontal axis, $H+(\lambda+\rho)(H),\;\;H\in \h$ on the vertical axis (as a function of degree 1), and $1$ at the corner of the diagram. By definition, $T_2(F)=1\ast_2 F$, and:
\begin{equation}\label{by}
(H+(\lambda+\rho)(H)) \ast _1 T_2(F) = (H+(\lambda+\rho)(H)) \ast _1 (1\ast _2 F)=((H+(\lambda+\rho)(H)) \ast _1 1) \ast _2 F = 0,
\end{equation}
by the compatibility relation (\ref{compatibility}) and Lemma \ref{lemma for H}. We prove first that $(H+(\lambda+\rho)(H))\ast_{DK} (\overline{T}_1^{-1}T_2(F))\in S_{(\epsilon)}(\g)\ast_{DK}\h_{\lambda+\rho}$. Recall that for $P,G \in (S(\g),\ast_{DK}),  K\in H^0_{(\epsilon)}(\h_{\lambda}^{\bot},d^{(\epsilon)}_{\g^{\ast},\h_{\lambda}^{\bot},\mathfrak{q}})$ we have $(P\ast_{DK} G) \ast _1 K = P\ast_1 (G\ast_1 K)$. Then, $\left((H+(\lambda+\rho)(H))\ast_{DK} (\overline{T}_1^{-1}T_2(F))\right)\ast_1 1 = (H+(\lambda+\rho)(H))\ast_1(\overline{T}_1^{-1}T_2(F)\ast_1 1)= \left(H+(\lambda+\rho)(H)\right)\ast_1 \left(\overline{T}_1(\overline{T}_1^{-1}T_2(F))\right) = (H+(\lambda+\rho)(H))\ast_1T_2(F)=0$,
where for the last equality we used (\ref{by}). 
Thus, $(H+(\lambda+\rho)(H))\ast_{DK} (\overline{T}_1^{-1}T_2(F))\in S_{(\epsilon)}(\g)\ast_{DK} \h_{\lambda+\rho}$ by Lemma \ref{lemma for T_1}. Using the decompositions (\ref{decomposition U(g)[[e]]}) and (\ref{decomposition S(g)[[e]]}) we get $\left((H+(\lambda+\rho)(H))\ast_{DK} \overline{T}_1^{-1}T_2(F) - \overline{T}_1^{-1}T_2(F)\ast_{DK} (H+(\lambda+\rho)(H))\right) \in S_{(\epsilon)}(\g)\ast_{DK} \h_{\lambda+\rho}$, that is, $ [(H+(\lambda+\rho)(H)),\overline{T}_1^{-1}T_2(F)]\in S_{(\epsilon)}(\g)\ast_{DK} \h_{\lambda+\rho}$. Since $[H+(\lambda+\rho)(H),\bullet]=[H,\bullet]$, we have that for $F\in H^0_{(\epsilon)}(\h_{\lambda}^{\bot},d^{(\epsilon)}_{\h_{\lambda}^{\bot},\mathfrak{q}})$, $\overline{T}_1^{-1}T_2(F)\in (S_{(\epsilon)}(\g)/S_{(\epsilon)}(\g)\ast_{DK} \h_{\lambda+\rho})^{\h},$ and $\overline{\beta}_{\mathfrak{q},(\epsilon)}\circ\partial_{q_{(\epsilon)}^{\frac{1}{2}}}\circ \overline{T}_1^{-1}T_2(F)\in(U_{(\epsilon)}(\g)/U_{(\epsilon)}(\g)\cdot \h_{\lambda+\rho})^{\h},$ by the beginning of $\mathcal{x}$ 5. $\diamond$
\newtheorem{okae2}[con]{Proposition}
\begin{okae2}\label{reverse main}
With these notations, $(U_{(\epsilon)}(\g)/U_{(\epsilon)}(\g)\h_{\lambda+\rho})^{\h} \subset H^0_{(\epsilon)}(\h^{\bot}_{\lambda},d^{(\epsilon)}_{\h^{\bot}_{\lambda},\mathfrak{q}})$.
\end{okae2}
\textit{Proof}. The proof of this fact is based on the description of the graphs that appear in a specific Stokes equation. This equation produces a system of homogeneous equations with respect to $\mathrm{deg}_{\epsilon}$, which will turn out to be the system  $d^{(\epsilon)}_{\h^{\bot}_{\lambda},\mathfrak{q}}(F)=0$ defining $ H^0_{(\epsilon)}(\h^{\bot}_{\lambda},d^{(\epsilon)}_{\h^{\bot}_{\lambda},\mathfrak{q}})$. Put again in the biquantization diagram, $\g^{\ast}$ on the vertical axis, and $-f+\h^{\bot}$ on the horizontal one. Let $G$ be an element of $(U_{(\epsilon)}(\g)/U_{(\epsilon)}(\g)\h_{\lambda+\rho})^{\h}$ and place it on the vertical axis. Let $H+(\lambda+\rho)(H),\;H\in \h$ be a function of degree 1 and place it also on the vertical axis. Since $(H+(\lambda+\rho)(H))\ast_{DK} G - G\ast_{DK} (H+(\lambda+\rho)(H)) \in S_{(\epsilon)}(\g)\ast_{DK} \h_{\lambda+\rho}$ by the proof of Proposition \ref{direct main}, we have $\left((H+(\lambda+\rho)(H))\ast_{DK} G - G\ast_{DK} (H+(\lambda+\rho)(H))\right)\ast_1 1 = 0$ by Lemma \ref{lemma for T_1}. Thus $((H+(\lambda+\rho)(H))\ast_{DK} G)\ast_1 1 =0,\; \textlatin{since} \;\;(H+(\lambda+\rho)(H))\ast_1 1 =0 \Rightarrow (G\ast_{DK} (H+(\lambda+\rho)(H)))\ast_1 1=0$. Finally, $((H+(\lambda+\rho)(H))\ast_{DK} G)\ast_1 1= (H+(\lambda+\rho)(H))\ast_1 (G\ast_1 1)= (H+(\lambda+\rho)(H))\ast_1[1\ast_2(T^{-1}_2\overline{T}_1(G))]=0$. Now set $T^{-1}_2\overline{T}_1(G)=F$ for an $F$ on the horizontal axis, and write

\begin{equation}\label{H1F=0}
(H+(\lambda+\rho)(H))\ast_1(1\ast_2 F)=0.
\end{equation} 

Let $s$ be a point on the horizontal axis. This point will be considered as $F$ and we will use it to calculate the following integral: Let $s$ move on the horizontal axis , and fix $1$ 
at the corner of the diagram. Let $C(s)= \sum_{\Gamma}\overline{\omega}_{\Gamma}(s)B_{\Gamma}(F)$ where $\Gamma$ are all the possible graphs appearing in the expression $(H+(\lambda+\rho)(H))\ast_1 1 \ast_2 F$.
To apply the Stokes equation, we need to calculate the integral $\sum _{\Gamma}\int^{\infty}_{0}\mathrm{d}\overline{\omega}_{\Gamma}(s)B_{\Gamma}(F)$ and for this, to calculate $\lim_{s\rightarrow \infty}C(s)$ and $\lim_{s\rightarrow 0}C(s)$ (where $0$ is considered to be the corner of the diagram).
\begin{figure}[h!]
\begin{center}
\includegraphics[width=6cm]{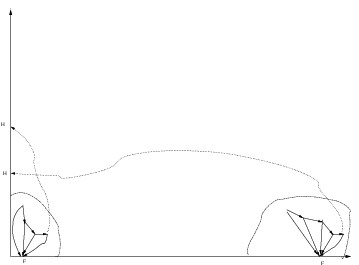}
\caption{\footnotesize Behaviour of $(H+(\lambda+\rho)(H))\ast_1 1 \ast_2 F$ when $s\rightarrow 0$ and $s\rightarrow \infty$ respectively.}
\end{center}
\end{figure}
First, when $s\longrightarrow 0$, the corresponding expression $(H+(\lambda+\rho)(H))\ast_1 (1\ast_2 F)$ is zero by (\ref{H1F=0}).
On the other side, when $s\longrightarrow \infty$ on the $\h_{\lambda}^{\bot}-$ axis,  is like $(H+(\lambda+\rho)(H))$ tending to the corner and the corresponding term $((H+(\lambda+\rho)(H))\ast_1 1)\ast_2 F$ is zero by Lemma \ref{lemma for H}. Recall that $H+(\lambda+\rho)(H)$ is of degree 1, so it can be derived by only one edge. The function $H+(\lambda+\rho)(H)$ is always derived by an edge since in the opposite case there will be an edge with no available endpoint. This will be clear after the description of the possible interior and exterior graphs that follows.

\textbf{The Stokes equation.} Let $F=\sum_{i=0}^n\epsilon^iF_{n-i}$ where the terms $F_j$ are not homogeneous polynomials. We saw before what happens when $s\rightarrow 0$ and $s\rightarrow \infty$ on the horizontal axis, so
\begin{equation}\label{first stokes}
\sum _{\Gamma}\int^{\infty}_{0}\mathrm{d}\overline{\omega}_{\Gamma}(s)B_{\Gamma}(F) = \sum_{\Gamma}(\overline{\omega}_{\Gamma}(\infty)B_{\Gamma}(F)- \overline{\omega}_{\Gamma}(0)B_{\Gamma}(F))=0-0=0,
\end{equation}
where $\Gamma$ runs over all the possible graphs. To write the differential $\mathrm{d}\overline{\omega}_{\Gamma}(s)$, we check the possible concentrations of points and graphs (both interior and exterior) so that the final contribution is nonzero. 
 
$\textbf{Possible concentrations appearing at the Stokes equation.}$ We first examine four kinds of concentrations that turn out to have zero contribution. These are the following:
\begin{enumerate}
\item \textbf{Concentrations at} $\mathbf{(H+(\lambda+\rho)(H))}.$ Concentrating at $(H+(\lambda+\rho)(H))$ at the vertical axis, we necessarily have a colored edge leaving the concentration. Indeed, since $(H+(\lambda+\rho)(H))$ is of degree 1 it can receive only one edge and the only possible graph in the concentration is the $\mathcal{B}$-type graph of one edge with another one $e_{\infty}$ leaving the concentration. This isn't possible as $e_{\infty}$ should be labeled here with a variable in $(\g^{\ast})^{\bot}$.

\item \textbf{Concentrations outside} $\mathbf{(H+(\lambda+\rho)(H))}$ \textbf{on the vertical axis.} Such concentrations leave two edges towards $\infty$. As before there are no variables for these edges and so the contribution is $0$.

\item \textbf{Aerial concentrations}. These may consist of concentrations of three or more vertices of type I, or of two such vertices collapsing. In this case we can regroup these terms in such way, so that their total contribution is zero as a result of the Jacobi identity.

\item \textbf{Concentrations away from s on the horizontal axis}. If we concentrate on the horizontal axis but not at $s$, the only possible graph will be the one of two type I vertices with two edges of the same color from one to the other and two colored edges leaving the concentration. However these terms have zero contribution due to Kontsevich's Lemma ~\cite{K}  $\mathcal{x}$ 7.3.3.1. 
\end{enumerate}

\begin{figure}[h!]
\begin{center}
\includegraphics[width=6cm]{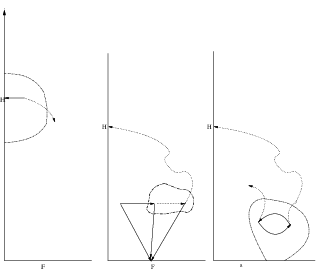}
\caption{\footnotesize The types 1,3 and 4 of concentrations that have zero contribution.}
\end{center}
\end{figure}

One can  have other, called here \textsl{exterior}, graphs acting on \textlatin{interior} concentrations. Let $\Gamma^{\alpha}_{int}$ be the family of possible interior graphs with $\alpha$ being the edge leaving the graph, and $\Gamma^{\alpha}_{ext}$ be the exterior ones acting on them. These graphs both depend on the edge $\alpha$. For graphs $\Gamma^{\alpha,k}_{int}\in \Gamma^{\alpha}_{int}$ and $\Gamma^{\alpha,s}_{ext}\in\Gamma^{\alpha}_{ext}$, the numbers $k,s \in \mathbb{Z}^+$ will denote the number of type I vertices in each case. This refers to the total number of type l vertices for interior and exterior graphs respectively.\footnote{For the exterior graphs, type II vertices are considered those representing the concentration so as type I vertices are considered all the others.}

\textbf{Interior graphs. (Denoted as} $\mathbf{\Gamma^{\alpha}_{int})}$. The principle for having a non-zero contribution by a concentration is that the dimension of the concentrated manifold $\overline{C}^+_{k,1}$ should be equal to the number of edges of a possible $\Gamma^{\alpha,k}_{int}$. So interior graphs can be only of $\mathcal{B}$-type or $\mathcal{BW}$-type. This is true by exclusion of the other possibility, $\mathcal{W}$-type graphs. Indeed, the concentration of a $\mathcal{W}$-type graph with $k$ vertices, will give an extra dimension to the concentration manifold, that is, any such graph will give a zero coefficient $\overline{\omega_{\Gamma}}$. In other words, after the concentration of $\Gamma^{\alpha,k}_{int}$ we won't have the right dimension for the interior concentation manifold. On the contrary, the two cases of graphs that we accept, get through this anomaly by allowing a colored edge $e_{\infty}$ to leave the concentration. 

\textbf{Exterior graphs. (Denoted as} $\mathbf{\Gamma^{\alpha}_{ext})}$
As a first possibility here we may have one $\mathcal{B}$-type graph receiving at its root the edge $\alpha=e_{\infty}$ from $\Gamma^{\alpha,k}_{int}$ and having a colored edge deriving $H+(\lambda+\rho)(H)$, plus an infinite number of superposing $\mathcal{W}$-type graphs deriving also the concentration. The other possibility is to have an arbitrary finite number of $\mathcal{W}$-type graphs. In this second case, the edge $\alpha$ leaving the interior concentration derives directly the function $H+(\lambda+\rho)(H)$. The case of a $\mathcal{BW}$-type graph is not acceptable as a possible $\Gamma^{\alpha}_{ext}$ for the following reason: The colored edge $e_{\infty}$ leaving a $\mathcal{BW}$-type graph will in that case derive directly the function $H+(\lambda+\rho)(H)$ on the vertical axis, so the edge $\alpha$ of $\Gamma^{\alpha,k}_{int}$ will have nothing to derive. As we saw before we have another possiblity for $\Gamma^{\alpha,k}_{int}$, namely having $\mathcal{W}$-type graphs. It can be shown (\cite{BAT} Lemma 3.7) that the existence of such graphs does not alter the equations (\ref{general}) so $\Gamma^{\alpha,k}_{int}$ will always be here a $\mathcal{B}-$ type graph with $k$ vertices of type I. For the exterior graphs now, we have also $\mathcal{B}$-type graphs in $\Gamma^{\alpha,s}_{ext}$, that receive at their root the colored edge coming from $\Gamma^{\alpha,k}_{int}$. So $\Gamma^{\alpha,m}_{ext}$ are $\mathcal{B}$-type graphs with $m$ type I vertices. The other possibility is to have additionally an arbitrary finite number of $\mathcal{W}$-type graphs.  If additionally to a $\mathcal{B}_m-$ type graph, we have $r$ $\mathcal{W}$-type graphs of $p_1,\ldots p_r$ type II vertices respectively, we will  denote this family of graphs as $\Gamma^{\alpha,m,p_1,\ldots,p_r}_{ext}$.

\begin{figure}[h!]
\begin{center}
\includegraphics[width=6cm]{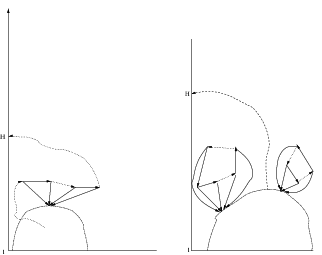}
\caption{\footnotesize The possible exterior graphs ($\Gamma^{\alpha}_{ext}$).}
\end{center}
\end{figure}

\begin{figure}[h!]
\begin{center}
\includegraphics[width=5cm]{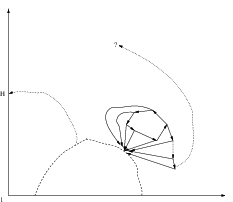}
\caption{\footnotesize $\mathcal{BW}-$ type graphs are excluded from the possible $\Gamma^{\alpha}_{ext}$.}
\end{center}
\end{figure}

\textbf{Writing the equations.}
We skip the notation $\overline{\omega}_{\Gamma}$ and just write $\omega_{\Gamma}$ since all graphs are colored. Let $\Gamma^{\alpha,k}_{int}\in\Gamma^{\alpha}_{int}$ and $\Gamma^{\alpha,m}_{ext}\in\Gamma^{\alpha}_{ext}$ and $B^k_{\Gamma^{\alpha}_{int}}$, $B^m_{\Gamma^{\alpha}_{ext}}$  be the respective differential operators. The equation we want to express in graphs is
\begin{equation}\label{Stokes}
\sum_{\Gamma^{\alpha}_{int},\Gamma^{\alpha}_{ext}}\int_0^{\infty}\hat{\omega}_{\Gamma^{\alpha}_{ext}}(s)B_{\Gamma^{\alpha}_{ext}}\left(\omega_{\Gamma^{\alpha}_{int}}B_{\Gamma^{\alpha}_{int}}\right)\mathrm{d}s=0
\end{equation}
The hat denotes the fact that $\hat{\omega}_{\Gamma^{\alpha}_{ext}}(s)$ is an 1- form while $\omega_{\Gamma^{\alpha}_{int}}$ is a coefficient. Note also that the term $\omega_{\Gamma^{\alpha}_{int}}B_{\Gamma^{\alpha}_{int}}$ doesn't depend on $s$. The idea is to analyze equation (\ref{Stokes}) in terms of the graphs $\Gamma^{\alpha}_{int}$ and $\Gamma^{\alpha}_{ext}$ for interior and exterior concentrations respectively.  As a convention, let $\Gamma^{\alpha,1}_{int}$ be the $\mathcal{B}$-type graph with one vertex and $\Gamma^{\alpha,0}_{ext}$ the $\mathcal{B}$-type graph with no vertices i.e an edge leaving from the concentration at the horizontal axis and ending at $H+(\lambda+\rho)(H)$.\\
\noindent We begin examining the possible concentrations of points and edges excluding the non-acceptable cases discussed before. The simplest concentration is to have the type I vertex of $\Gamma^{\alpha,1}_{int}$ collapsing to the horizontal axis, and to have no other type I vertices for the exterior graph. The resulting graph is $\Gamma^{\alpha,0}_{ext}$, the situation before the collapse is depicted in Fig. \ref{graph degree 1}.  So the first term in the left hand sum of equation (\ref{Stokes}) is $\int_0^\infty\hat{\omega}_{\Gamma^{\alpha,0}_{ext}}(s)B_{\Gamma^{\alpha,0}_{ext}}(\omega_{\Gamma^{\alpha,1}_{int}}
B_{\Gamma^{\alpha,1}_{int}})\mathrm{d}s$. It can be computed (see ~\cite{BAT} Lemma 3.7) that $B_{\Gamma^{\alpha,0}_{ext}}$ is a nonzero number and in particular $B_{\Gamma^{\alpha,0}_{ext}}=
\omega_{\Gamma^{\alpha,0}_{ext}}=1$. Thus this term equals $\omega_{\Gamma^{\alpha,1}_{int}}
B_{\Gamma^{\alpha,1}_{int}}$ where $\omega_{\Gamma^{\alpha,1}_{int}}=\int_{\overline{C}^+_{1,1}}
\mathrm{d}\phi_{+,e}$. When we have in total two type II vertices, then one has the following possible concentrations: First, to have as the interior graph the graph $\Gamma^{\alpha,1}_{int}$ collapsing to the horizontal axis and as exterior, the graph $\Gamma^{\alpha,1}_{ext}$. The second possible case is to collapse the graph $\Gamma^{\alpha,2}_{int}$ at the horizontal axis, and have again $\Gamma^{\alpha,0}_{ext}$ as the exterior graph.  These situations give the summand  
\begin{equation}\label{second summand}
\int_0^\infty\hat{\omega}_{\Gamma^{\alpha,1}_{ext}}(s)B_{\Gamma^{\alpha,1}_{ext}}
(\omega_{\Gamma^{\alpha,1}_{int}}
B_{\Gamma^{\alpha,1}_{int}})\mathrm{d}s+
\int_0^\infty\hat{\omega}_{\Gamma^{\alpha,0}_{ext}}(s)B_{\Gamma^{\alpha,0}_{ext}}
(\omega_{\Gamma^{\alpha,2}_{int}}
B_{\Gamma^{\alpha,2}_{int}})\mathrm{d}s.
\end{equation}
at the left hand side of (\ref{Stokes}). Following the previous calculation, this term is $ \omega_{\Gamma^{\alpha,2}_{int}}
B_{\Gamma^{\alpha,2}_{int}}+\omega_{\Gamma^{\alpha,1}_{ext}}
B_{\Gamma^{\alpha,1}_{ext}}(\omega_{\Gamma^{\alpha,1}_{int}}
B_{\Gamma^{\alpha,1}_{int}})$, the new coefficient being $\omega_{\Gamma^{\alpha,2}_{int}}=\int_{\overline{C}^+_{2,1}}
\mathrm{d}\phi_{+,e_1^1}\wedge
\mathrm{d}\phi_{-,e_1^2}\wedge\mathrm{d}\phi_{+,e_2^2}$, if we put the root of the $\mathcal{B}_2$ graph in first position, the edge leaving the root and deriving the axis in first position, and the edge $e_{\infty}$ in first position with respect to the edges leaving the second vertex. Computing the possible concentrations for larger numbers of type I vertices, one gets inductively the sum at the left hand side of (\ref{Stokes}). This sum can be grouped in nice terms if we take in mind the parameter $\epsilon$ for each type I vertex. In fact one can write homogeneous equations with respect to the total $\deg_{\epsilon}$ degree. If the total graph derives the function $F=\sum_i^n \epsilon^iF_{n-i}$, in differential operator terms, the equation (\ref{Stokes}) is equivalent to

\begin{equation}\label{general}
\sum_{\alpha}\left(\sum_{\Gamma^{\alpha}_{int},\Gamma^{\alpha}_{ext}} \left(\omega_{\Gamma^{\alpha}_{ext}}B_{\Gamma^{\alpha}_{ext}}(\omega_{\Gamma^{\alpha}_{int}}B_{\Gamma^{\alpha}_{int}}(F))\right)\right)=0 
\end{equation}
\[\Leftrightarrow \sum_{\alpha}\left(\sum_{\Gamma^{\alpha}_{int},\Gamma^{\alpha}_{ext}}\left(\sum_{l,k,m}
\left(\omega_{\Gamma^{\alpha,m}_{ext}}B_{\Gamma^{\alpha,m}_{ext}}(\omega_{\Gamma^{\alpha,k}_{int}}B_{\Gamma^{\alpha,k}_{int}}(F_l))\right)\epsilon^{m+k+l}\right)\right)=0,\;\;\;\textlatin{for}\;\; l=1,\ldots n, k,m=0,\ldots \infty\]

\begin{figure}[h!]\label{graph 2 first}
\begin{center}
\includegraphics[width=8cm]{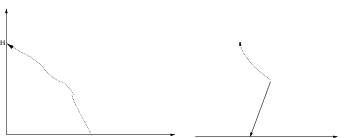}
\caption{\footnotesize The operators $B_{\Gamma^{\alpha,0}_{ext}}$ and $B_{\Gamma^{\alpha,1}_{int}}$.}
\end{center}
\end{figure}
\noindent For $\mathrm{deg}_{\epsilon}=1$, the only term is $\omega_{\Gamma^{\alpha,0}_{ext}}B_{\Gamma^{\alpha,0}_{ext}}(\omega_{\Gamma^{\alpha,1}_{int}}B_{\Gamma^{\alpha,1}_{int}})$ as it is clear from Fig. \ref{graph 2 first}. Thus, our first equation is $\omega_{\Gamma^{\alpha,0}_{ext}}B_{\Gamma^{\alpha,0}_{ext}}(\omega_{\Gamma^{\alpha,1}_{int}}B_{\Gamma^{\alpha,1}_{int}}(F_n))=0$, 
\begin{figure}[h!]\label{graph degree 1}
\begin{center}
\includegraphics[width=4cm]{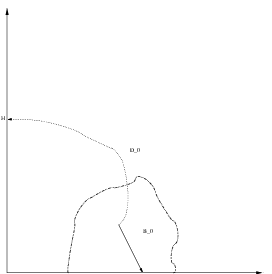}
\caption{\footnotesize The graph corresponding to $B^0_{\Gamma^{\alpha}_{ext}}(B^1_{\Gamma^{\alpha}_{int}})$.}\label{D_0(B_0)}
\end{center}
\end{figure}
and $\omega_{\Gamma^{\alpha,0}_{ext}}B_{\Gamma^{\alpha,0}_{ext}}(\omega_{\Gamma^{\alpha,1}_{int}}B_{\Gamma^{\alpha,1}_{int}}(\cdot))$ is an operator of degree 1.  So this first equation yields
\begin{equation}\label{degree 1}
\omega_{\Gamma^{\alpha,1}_{int}}B_{\Gamma^{\alpha,1}_{int}}(F_n)=0.
\end{equation}
Taking in mind the graphs that we allow for 
$\Gamma^{\alpha,1}_{int}$, this is actually the first equation $d^{(1)}_{\h_{\lambda}^{\bot},\mathfrak{q}}(F_n)=0$ of the system $d^{(\epsilon)}_{\h^{\bot}_{\lambda},\mathfrak{q}}(F)=0$ defining the reduction algebra $H^0_{(\epsilon)}(\h^{\bot}_{\lambda},d^{(\epsilon)}_{\h^{\bot}_{\lambda},\mathfrak{q}})$. For $\mathrm{deg}_{\epsilon}=2$ we have $\sum_{\alpha}\left(\omega_{\Gamma^{\alpha,1}_{ext}}B_{\Gamma^{\alpha,1}_{ext}}(\omega_{\Gamma^{\alpha,1}_{int}}B_{\Gamma^{\alpha,1}_{int}}(F_n))+\omega_{\Gamma^{\alpha,0}_{ext}}B_{\Gamma^{\alpha,0}_{ext}}(\omega_{\Gamma^{\alpha,2}_{int}}B_{\Gamma^{\alpha,2}_{int}}(F_n)) + \omega_{\Gamma^{\alpha,0}_{ext}}B_{\Gamma^{\alpha,0}_{ext}}(\omega_{\Gamma^{\alpha,1}_{int}}B_{\Gamma^{\alpha,1}_{int}}(F_{n-1}))\right)=0$. By (\ref{degree 1}) we get $\left(\omega_{\Gamma^{\alpha,2}_{int}}B_{\Gamma^{\alpha,2}_{int}}(F_n) + \omega_{\Gamma^{\alpha,1}_{int}}B_{\Gamma^{\alpha,1}_{int}}(F_{n-1})\right)\omega_{\Gamma^{\alpha,0}_{ext}}\partial_{\alpha}H = 0 \;\;\; \forall \alpha, \forall H$ which gives that $\forall\alpha$,

\begin{equation}\label{degree 2}
\omega_{\Gamma^{\alpha,2}_{int}}B_{\Gamma^{\alpha,2}_{int}}(F_n) +\omega_{\Gamma^{\alpha,1}_{int}} B_{\Gamma^{\alpha,1}_{int}}(F_{n-1})=0\;.
\end{equation}

\noindent Again, we observe that this corresponds to the second equation $d^{(2)}_{\h_{\lambda}^{\bot},\mathfrak{q}}(F_n)+d^{(1)}_{\h_{\lambda}^{\bot},\mathfrak{q}}(F_{n-1})=0$ of the system $d^{(\epsilon)}_{\h^{\bot}_{\lambda},\mathfrak{q}}(F)=0$. For $\mathrm{deg}_{\epsilon}=3$, we write 
\begin{eqnarray}\label{general degree 3}
\sum_{\alpha}[\omega_{\Gamma^{\alpha,0}_{ext}}B_{\Gamma^{\alpha,0}_{ext}}(\omega_{\Gamma^{\alpha,3}_{int}}B_{\Gamma^{\alpha.3}_{int}}(F_n)) + \omega_{\Gamma^{\alpha,2}_{ext}}B_{\Gamma^{\alpha,2}_{ext}}(\omega_{\Gamma^{\alpha,1}_{int}}B_{\Gamma^{\alpha,1}_{int}}(F_n)) +\omega_{\Gamma^{\alpha,1}_{ext}} B_{\Gamma^{\alpha,1}_{ext}}(\omega_{\Gamma^{\alpha,2}_{int}}B_{\Gamma^{\alpha,2}_{int}}(F_n)) + 
\nonumber \\
+\omega_{\Gamma^{\alpha,1}_{ext}}B_{\Gamma^{\alpha,1}_{ext}}(\omega_{\Gamma^{\alpha,1}_{int}}B_{\Gamma^{\alpha,1}_{int}}(F_{n-1})) + \omega_{\Gamma^{\alpha,0}_{ext}}B_{\Gamma^{\alpha,0}_{ext}}(\omega_{\Gamma^{\alpha,2}_{int}}B_{\Gamma^{\alpha,2}_{int}}(F_{n-1})) +\omega_{\Gamma^{\alpha,0}_{ext}} B_{\Gamma^{\alpha,0}_{ext}}(\omega_{\Gamma^{\alpha,1}_{int}}B_{\Gamma^{\alpha,1}_{int}}(F_{n-2}))] = 0
\end{eqnarray}

\noindent Now by (\ref{degree 2}) we get $[\omega_{\Gamma^{\alpha,3}_{int}}B_{\Gamma^{\alpha,3}_{int}}(F_n) +\omega_{\Gamma^{\alpha,2}_{int}} B_{\Gamma^{\alpha,2}_{int}}(F_{n-1}) + \omega_{\Gamma^{\alpha,1}_{int}}B_{\Gamma^{\alpha,1}_{int}}(F_{n-2})]\omega_{\Gamma^{\alpha,0}_{ext}}\partial_{\alpha}H=0 \;\;\; \forall \alpha,\forall H$ and as above, for all indices $\alpha$ this gives
\begin{equation}\label{degree 3}
\omega_{\Gamma^{\alpha,3}_{int}}B_{\Gamma^{\alpha,3}_{int}}(F_n) +\omega_{\Gamma^{\alpha,2}_{int}} B_{\Gamma^{\alpha,2}_{int}}(F_{n-1}) + \omega_{\Gamma^{\alpha,1}_{int}}B_{\Gamma^{\alpha,1}_{int}}(F_{n-2})=0.
\end{equation}
Again, considering the graphs that we allow for $\Gamma^{\alpha}_{int}$, this last equation is actually the third equation $d^{(3)}_{\h_{\lambda}^{\bot},\mathfrak{q}}(F_n)+d^{(2)}_{\h_{\lambda}^{\bot},\mathfrak{q}}(F_{n-1})+d^{(1)}_{\h_{\lambda}^{\bot},\mathfrak{q}}(F_{n-2})=0$ of the system $d^{(\epsilon)}_{\h^{\bot}_{\lambda},\mathfrak{q}}(F)=0$. Induction on $\mathrm{deg}_{\epsilon}$ will give, as (\ref{degree 1}), (\ref{degree 2}), (\ref{degree 3})  suggest, 

\[\sum_{\alpha}\omega_{\Gamma^{\alpha}_{ext}}B_{\Gamma^{\alpha}_{ext}}(\omega_{\Gamma^{\alpha}_{int}}B_{\Gamma^{\alpha}_{int}}(F))=0 \Leftrightarrow \sum_{\alpha}\sum_{l,i,j}\omega_{\Gamma^{\alpha,j}_{ext}}B_{\Gamma^{\alpha,j}_{ext}}(\omega_{\Gamma^{\alpha,i}_{int}}B_{\Gamma^{\alpha,i}_{int}}(F_l))=0 \Leftrightarrow \sum_i\epsilon^id^{(i)}_{\h_{\lambda}^{\bot},\mathfrak{q}}(F)=0\Leftrightarrow d^{(\epsilon)}_{\h_{\lambda}^{\bot},\mathfrak{q}}(F)=0.\]

\noindent Thus for $F\in (U_{(\epsilon)}(\g)/U_{(\epsilon)}(\g)\h_{\lambda+\rho})^{\h}$ the Stokes equation (\ref{first stokes}) produces the reduction equations for $H^0_{(\epsilon)}(\h^{\bot}_{\lambda},d^{(\epsilon)}_{\h^{\bot}_{\lambda},\mathfrak{q}})$, that is $G=\left(\overline{\beta}_{\mathfrak{q},(\epsilon)}\circ\partial_{q_{(\epsilon)}^{\frac{1}{2}}}\circ \overline{T}_1^{-1}T_2\right)(F) \in (U_{(\epsilon)}(\g)/U_{(\epsilon)}(\g)\h_{\lambda+\rho})^{\h} \Rightarrow F\in H^0_{(\epsilon)}(\h^{\bot}_{\lambda},d^{(\epsilon)}_{\h^{\bot}_{\lambda},\mathfrak{q}})$.

\section{Specialization and deformations.}
\textbf{6.1 Homogeinity degree and the vector space case.} At the proof of theorem \ref{main} the deformation parameter $\epsilon$ was used to write down homogeneous equations. These equations were a product of a Stokes equation (\ref{Stokes}), and resulted in the reduction equations of $H^0_{(\epsilon)}(\h_{\lambda}^{\bot},d^{(\epsilon)}_{\h_{\lambda}^{\bot},\mathfrak{q}})$. Recall that in $\mathcal{x}$ 3 we used the parameter $\epsilon$ to describe the differential $d^{(\epsilon)}_{\h^{\bot},\mathfrak{q}}$ and give the definition of $H^0_{(\epsilon)}(\h^{\bot},d^{(\epsilon)}_{\h^{\bot},\mathfrak{q}})$. However it was explained later that in the vector space case (i.e $\h^{\bot}$), the algebra $H^0(\h^{\bot},d_{\h^{\bot},\mathfrak{q}})$ could be defined without using $\epsilon$, but using instead the degree $\mathrm{deg}_{\mathfrak{q}}$ of operators in $d_{\h^{\bot},\mathfrak{q}}$. 
Thus for a Lie algebra $\g$, a subalgebra $\h\subset\g$ and $\mathfrak{q}$ such that $\g=\h\oplus\mathfrak{q}$, one has
\[(U(\g)/U(\g)\h)^{\h}\simeq H^0(\h^{\bot},d_{\h^{\bot},\mathfrak{q}})\simeq H^0_{(\epsilon=1)} (\h^{\bot},d^{(\epsilon=1)}_{\h^{\bot},\mathfrak{q}}).\]

The second isomorphism is Lemma \ref{ftou}. The direction $H^0(\h^{\bot},d_{\h^{\bot},\mathfrak{q}})\subset (U(\g)/U(\g)\h)^{\h}$ can be proved similarly to Theorem \ref{main} setting $\lambda=0$. For the opposite direction of the first isomorphism, $(U(\g)/U(\g)\h)^{\h}\subset H^0(\h^{\bot},d_{\h^{\bot},\mathfrak{q}})$, we have again interior and exterior graphs describing the Stokes equation (\ref{Stokes}). This time we argue in terms of the degree $\mathrm{deg}_{\mathfrak{q}}$ of the operators $B_{\Gamma^{\alpha}_{ext}},\;B_{\Gamma^{\alpha}_{ext}}$. Namely we use the total degree in $\mathfrak{q}$ to write down the equations (\ref{general}) (as we used the total degree in $\epsilon$ in the proof of Theorem \ref{main}). For $\mathrm{deg}_{\mathfrak{q}}(B_{\Gamma^{\alpha}_{ext}})+\mathrm{deg}_{\mathfrak{q}}(B_{\Gamma^{\alpha}_{int}})=1$ we get $B^0_{\Gamma^{\alpha}_{ext}}(B^1_{\Gamma^{\alpha}_{int}}(F_n))=0,$ and then $B^1_{\Gamma^{\alpha}_{int}}(F_n)=0$. Similarly we get the rest of the equations in (\ref{general}) for $\mathrm{deg}_{\mathfrak{q}}(B_{\Gamma^{\alpha}_{ext}})+\mathrm{deg}_{\mathfrak{q}}(B_{\Gamma^{\alpha}_{int}})=k$. 

Keeping  the notation $\overline{T}_1:=T_1|_{S(\mathfrak{q})}$ also in the case without the parameter $\epsilon$, the map $\overline{\beta}_{\mathfrak{q}}\circ\partial_{q^{\frac{1}{2}}}\circ \overline{T}_1^{-1}T_2:\;H^0(\h_{\lambda}^{\bot},d_{\h_{\lambda}^{\bot},\mathfrak{q}})\hookrightarrow (U(\g)/U(\g)\h_{\lambda+\rho})^{\h},\;\;P\mapsto \overline{\beta}_{\mathfrak{q}}\circ\partial_{q^{\frac{1}{2}}}\circ \overline{T}_1^{-1}T_2(P)$,
is an injective algebra map since one simply has to work as in the first part of Theorem \ref{main}.

\textbf{6.2 Specialization algebra.} Recall the following notation from earlier: Let $\pi_{(\epsilon=1)}:\;H^0_{(\epsilon)}(\h_{\lambda}^{\bot},d^{(\epsilon)}_{\h^{\bot}_{\lambda},\mathfrak{q}})\longrightarrow H^0_{(\epsilon=1)}(\h_{\lambda}^{\bot},d^{(\epsilon=1)}_{\h^{\bot}_{\lambda},\mathfrak{q}})$ be the canonical projection. For $F^{'}=\sum_r\epsilon^rF^{'}_r,$ with $F^{'}\in  H^0_{(\epsilon)}(\h_{\lambda}^{\bot},d^{(\epsilon)}_{\h^{\bot}_{\lambda},\mathfrak{q}})$, we  denoted as $J:\;H^0_{(\epsilon)}(\h_{\lambda}^{\bot},d^{(\epsilon)}_{\h^{\bot}_{\lambda},\mathfrak{q}})\longrightarrow H^0(\h_{\lambda}^{\bot},d_{\h^{\bot}_{\lambda},\mathfrak{q}})$ the linear map defined by the formula $J(F^{'})=\sum_k F^{'}_k$.
Let also $\overline{J}:\;H^0_{(\epsilon=1)}(\h_{\lambda}^{\bot},d^{(\epsilon=1)}_{\h^{\bot}_{\lambda},\mathfrak{q}})\longrightarrow H^0(\h_{\lambda}^{\bot},d_{\h^{\bot}_{\lambda},\mathfrak{q}})$ be the quotient map (Lemma \ref{jbar}). In the following we will use the reduction algebra $H^0(\h_{t\lambda}^{\bot},d_{\h_{t\lambda}^{\bot},\mathfrak{q}})$, $t\in\mathbb{R}^{\ast}$. The operators in the differential defining $H^0(\h_{t\lambda}^{\bot},d_{\h_{t\lambda}^{\bot},\mathfrak{q}})$ are the same with the differential defining $H^0(\h_{\lambda}^{\bot},d_{\h_{\lambda}^{\bot},\mathfrak{q}})$. However in the case of $H^0(\h_{t\lambda}^{\bot},d_{\h_{t\lambda}^{\bot},\mathfrak{q}})$, the variable $t$ shows up after restriction of an operator at $-t\lambda+\h^{\bot}$ (this is due to the root of a $\mathcal{B}-$ type graph in $d_{\h_{t\lambda}^{\bot},\mathfrak{q}})$.

\newtheorem{okrat}[con]{Theorem}
\begin{okrat}\label{difo}
Let $\g$ be a Lie algebra, $\h\subset\g$ a subalgebra and $\mathfrak{q}$ such that $\g=\h\oplus\mathfrak{q}$. 
\begin{enumerate}
\item
 Let $F\in H^0(\h_{\lambda}^{\bot},d_{\h_{\lambda}^{\bot},\mathfrak{q}})$. Suppose an $F_{(t)}=\sum_pt^pF_p$ satisfying $F_{(t)}\in H^0(\h_{t\lambda}^{\bot},d_{\h_{t\lambda}^{\bot},\mathfrak{q}})$, $\forall t\in\mathbb{R}^{\ast}$ and $F_{(t=1)}=F$. Let $F_p=\sum_iF_p^{(i)}$ be a decomposition in homogeneous polynomials i.e $\mathrm{deg}_{\mathfrak{q}}(F^{(i)}_p)=i$ and let $F_{(\epsilon)}:=\epsilon^N\sum F_p^{(i)}\frac{1}{\epsilon^{i+p}}$ (with $N>>max(i+p)$). Then $F_{(\epsilon)}\in H^0_{(\epsilon)}(\h_{\lambda}^{\bot},d^{(\epsilon)}_{\h_{\lambda}^{\bot},\mathfrak{q}}),\;J(F_{(\epsilon)})=F$.
\item
Let $F\in H^0(\h_{\lambda}^{\bot},d_{\h_{\lambda}^{\bot},\mathfrak{q}})$. Suppose an $F_{(\epsilon)}=\sum_{0 \geq k\geq n}\epsilon^kF_k$ satisfying $F_{(\epsilon)}\in H^0_{(\epsilon)}(\h_{\lambda}^{\bot},d^{(\epsilon)}_{\h_{\lambda}^{\bot},\mathfrak{q}})$ and $J(F_{(\epsilon)})=F$. Let $F_k=\sum_iF_k^{(i)}$ be a decomposition in homogeneous polynomials i.e $\mathrm{deg}_{\mathfrak{q}}(F^{(i)}_k)=i$. Let also $F_{(t)}:=t^N\sum_{i,k}\frac{1}{t^{i+k}}F_k^{(i)}$ (with $N>>max(i+k)$). Then $\forall t\in\mathbb{R}^{\ast},\;\;F_{(t)}\in H^0(\h_{t\lambda}^{\bot},d_{\h_{t\lambda}^{\bot},\mathfrak{q}})$.
\end{enumerate}
\end{okrat}

\textbf{Proof of 1.}
Let's first examine the equations defining $H^0(\h_{t\lambda}^{\bot},d_{\h_{t\lambda}^{\bot},\mathfrak{q}})$. Let $B_k$ be an operator in $d^k_{\h_{t\lambda}^{\bot},\mathfrak{q}}$. Then $B_k|_{-t\lambda+\h^{\bot}}=tB_k'+B_k''$. Here $B_k'$ are operators whose graph is a $\mathcal{B}-$ type graph with $k$ type I vertices (i.e of $\mathcal{B}_k-$ type) and with its root in $\h^{\ast}$, while $B_k''$ are operators whose graph is either a $\mathcal{BW}_k-$ type graph, or a $\mathcal{B}_k-$ type graph with its root in $\mathfrak{q}^{\ast}$. Operators in $B_k'$ are of degree $-k$ (with respect to $\mathfrak{q}$ since they derive $k$ times). For operators in $B_k''$ we have that by definition $\mathcal{BW}_k-$ type graphs correspond to operators of degree $-k+1$. The same is true for graphs of $\mathcal{B}_k-$ type with their root in $\mathfrak{q}$.
Thus, $\mathrm{deg}_{\mathfrak{q}}(B_k')=-k$ and $\mathrm{deg}_{\mathfrak{q}}(B_k'')=-k+1$. Taking into consideration also the variable $t$, we can write that the total degree $\mathrm{deg}_{\mathfrak{q},t}:=\mathrm{deg}_{\mathfrak{q}}+\mathrm{deg}_{t}$ of $tB_k',B_k''$ is $-k+1$, that is $\mathrm{deg}_{\mathfrak{q},t}(tB_k')=\mathrm{deg}_{\mathfrak{q},t}(B_k'')=-k+1$. We decompose now the differential $d_{\h_{t\lambda}^{\bot},\mathfrak{q}}$ with respect to the total degree $\mathrm{deg}_{\mathfrak{q},t}$ of the operators in it. Namely we write $d_{\h_{t\lambda}^{\bot},\mathfrak{q}}=\sum_{k\geq 1}d^{|k|}_{\h_{t\lambda}^{\bot},\mathfrak{q}}$ where $\mathrm{deg}_{\mathfrak{q},t}(d^{|k|}_{\h_{t\lambda}^{\bot},\mathfrak{q}})=-k+1$. This sum can be decomposed as 

\begin{equation}\label{red tl}
d_{\h_{t\lambda}^{\bot},\mathfrak{q}}=\sum_{k\geq 1}td^{|k|^{'}}_{\h_{\lambda}^{\bot},\mathfrak{q}}+d^{|k|^{''}}_{\h_{\lambda}^{\bot},\mathfrak{q}}
\end{equation}
where $\mathrm{deg}_{\mathfrak{q}}(d^{|k|^{'}}_{\h_{\lambda}^{\bot},\mathfrak{q}})=-k$ and $\mathrm{deg}_{\mathfrak{q}}(d^{|k|^{''}}_{\h_{\lambda}^{\bot},\mathfrak{q}})=-k+1$. Let $G_{(t)}=\sum_pt^pF_p$ with $G_{(t)}\in  H^0(\h_{t\lambda}^{\bot},d_{\h_{t\lambda}^{\bot},\mathfrak{q}})$. Then $d_{\h_{t\lambda}^{\bot},\mathfrak{q}}(G_{(t)})=0$. So we have

\[d_{\h_{t\lambda}^{\bot},\mathfrak{q}}(G_{(t)})= 0\;\Leftrightarrow \sum_{k\geq 1} \left(td^{|k|^{'}}_{\h_{\lambda}^{\bot},\mathfrak{q}}+d^{|k|^{''}}_{\h_{\lambda}^{\bot},\mathfrak{q}}\right)\left(\sum_pt^pF_p\right)=0\Leftrightarrow \sum_pt^{p+1}\left(\sum_{k\geq 1} d^{|k|^{'}}_{\h_{\lambda}^{\bot},\mathfrak{q}}(F_p)+d^{|k|^{''}}_{\h_{\lambda}^{\bot},\mathfrak{q}}(F_{p+1})\right)=0\]

which is equivalent $\forall p\geq0,\;b\in \mathbb{Z}+$ to 

\begin{equation}\label{firsteq1}
\sum_{\substack{i,k\\i-k=b}}d^{|k|^{'}}_{\h_{\lambda}^{\bot},\mathfrak{q}}(F_p^{(i)})+d^{|k|^{''}}_{\h_{\lambda}^{\bot},\mathfrak{q}}(F_{p+1}^{(i-1)})=0.
\end{equation}
Set now $G_{(\epsilon)}=\sum F_p^{(i)}\frac{1}{\epsilon^{i+k}}$ with $\mathrm{deg}_{\mathfrak{q}}(F_p^{(i)})=i$. We want to calculate $d^{(\epsilon)}_{\h_{\lambda}^{\bot},\mathfrak{q}}(G_{(\epsilon)})$. This differential is written again grouping the family of $\mathcal{B}-$ type graphs with root in $\h^{\ast}$ and $\mathcal{BW}-$ type graphs together with $\mathcal{B}-$ type graphs having their root in $\mathfrak{q}^{\ast}$. Thus $d^{(\epsilon)}_{\h_{\lambda}^{\bot},\mathfrak{q}}$ can be written as

\begin{equation}\label{dife}
d^{(\epsilon)}_{\h_{\lambda}^{\bot},\mathfrak{q}}=\sum_{k\geq 1}\epsilon^kd^{(k)}_{\h_{\lambda}^{\bot},\mathfrak{q}}=\sum_{k\geq 1}\epsilon^k\left(d^{|k|^{'}}_{\h_{\lambda}^{\bot},\mathfrak{q}}+d^{|k|^{''}}_{\h_{\lambda}^{\bot},\mathfrak{q}}\right)
\end{equation}
 where $\mathrm{deg}_{\mathfrak{q}}(d^{|k|^{'}}_{\h_{\lambda}^{\bot},\mathfrak{q}})=-k$ and $\mathrm{deg}_{\mathfrak{q}}(d^{|k|^{''}}_{\h_{\lambda}^{\bot},\mathfrak{q}})=-k+1$. Denoting again as $B_k$  an operator in $d^{(k)}_{\h_{\lambda}^{\bot},\mathfrak{q}}$ we have

\[d^{(\epsilon)}_{\h_{\lambda}^{\bot},\mathfrak{q}}(G_{(\epsilon)})= \left(\sum_k \epsilon^k(\sum_{B_k\in d^{(k)}_{\h_{\lambda}^{\bot},\mathfrak{q}}}B_k)\right)\left(\sum_p F_p^{(i)}\frac{1}{\epsilon^{i+p}}\right)=\left(\sum_k\epsilon^k (d^{|k|^{'}}_{\h_{\lambda}^{\bot},\mathfrak{q}}+d^{|k|^{''}}_{\h_{\lambda}^{\bot},\mathfrak{q}})\right)\left(\sum_p F_p^{(i)}\frac{1}{\epsilon^{i+p}}\right)=\]
\[=\sum_{a}\epsilon^{a}\left(\sum_{b}\left(\sum_{\substack{i,k,p\\k-i-p=a\\i-k=b}}d^{|k|^{'}}_{\h_{\lambda}^{\bot},\mathfrak{q}}(F_p^{(i)})+d^{|k|^{''}}_{\h_{\lambda}^{\bot},\mathfrak{q}}(F_{p+1}^{(i-1)})\right)\right)\]

 for $a\in \mathbb{Z},b\in \mathbb{Z}^+$. For the summation index, we get $-p=a+b$, so the inner sum  for $p=-a-b$ is,
\begin{equation}\label{second eq}
\sum_{\substack{i,k\\i-k=b}}d^{|k|^{'}}_{\h_{\lambda}^{\bot},\mathfrak{q}}(F_p^{(i)})+d^{|k|^{''}}_{\h_{\lambda}^{\bot},\mathfrak{q}}(F_{p+1}^{(i-1)}).
\end{equation}
By (\ref{firsteq1}), the sum (\ref{second eq}) is $0$, and so $G_{(\epsilon)}\in  H^0_{(\epsilon)}(\h_{\lambda}^{\bot},d^{(\epsilon)}_{\h_{\lambda}^{\bot},\mathfrak{q}})$. Finally $F=F_{(t=1)}=\sum_{i,p}F_p^{(i)}=J(F_{(\epsilon)})$.

\textbf{Proof of 2.} We keep the notation as before. Let $F\in H^0(\h_{\lambda}^{\bot},d_{\h_{\lambda}^{\bot},\mathfrak{q}})$ and $F_{(\epsilon)}\in H^0_{(\epsilon)}(\h_{\lambda}^{\bot},d^{(\epsilon)}_{\h_{\lambda}^{\bot},\mathfrak{q}})$ with $F_{(\epsilon)}=F_0+\epsilon F_1+\epsilon^2 F_2+\cdots +\epsilon^nF_n$ and $F=J(F_{(\epsilon)})$.  For $k\in \{0,\ldots,n\}$ write $F_k$ as $F_k=\sum_iF_k^{(i)}$ where the $F_k^{(i)}$ are homogeneous polynomials of degree $\mathrm{deg}_{\mathfrak{q}}(F_k^{(i)})=i$. We have that 
\[d^{(\epsilon)}_{\h_{\lambda}^{\bot},\mathfrak{q}}(F_{(\epsilon)})=0\Leftrightarrow \left(\sum_k \epsilon^k(\sum_{B_k\in d^{(k)}_{\h_{\lambda}^{\bot},\mathfrak{q}}}B_k)\right)\left(\sum_p \epsilon^pF_p\right)=0\Leftrightarrow \left(\sum_k\epsilon^k (d^{|k|^{'}}_{\h_{\lambda}^{\bot},\mathfrak{q}}+d^{|k|^{''}}_{\h_{\lambda}^{\bot},\mathfrak{q}})\right)\left(\sum_{i,p}\epsilon^p F_p^{(i)}\right)=0\Leftrightarrow\]
\[\Leftrightarrow \sum_a\epsilon^a\left(\sum_b\left(\sum_{\substack{i,k\\k+p=a\\i-k=b}}d^{|k|^{'}}_{\h_{\lambda}^{\bot},\mathfrak{q}}(F_p^{(i)})+\sum_{\substack{i,k\\k+p=a\\i-k=b}}d^{|k|^{''}}_{\h_{\lambda}^{\bot},\mathfrak{q}}(F_p^{(i-1)})\right)\right)=0\]

which $\forall (a,b)\in\mathbb{Z}\times\mathbb{Z}^+$ is equivalent to

\begin{equation}\label{firsteq2}
\sum_{\substack{i,k\\k+p=a\\i-k=b}}d^{|k|^{'}}_{\h_{\lambda}^{\bot},\mathfrak{q}}(F_p^{(i)})+\sum_{\substack{i,k\\k+p=a\\i-k=b}}d^{|k|^{''}}_{\h_{\lambda}^{\bot},\mathfrak{q}}(F_p^{(i-1)})=0.
\end{equation}

Let now $G_t=\sum_{i,p}\frac{1}{t^{i+p}}F_p^{(i)}$. Then 
\[d_{\h_{t\lambda}^{\bot},\mathfrak{q}}(G_t)= \left(\sum_{k\geq 1}td^{|k|^{'}}_{\h_{\lambda}^{\bot},\mathfrak{q}}+d^{|k|^{''}}_{\h_{\lambda}^{\bot},\mathfrak{q}}\right)\left(\sum_{i,p}\frac{1}{t^{i+p}}F_p^{(i)}\right)=\sum_{\gamma}t^{-\gamma+1}\left(\sum_b\left(\sum_{\substack{i,k\\i+p=\gamma\\i-k=b}} d^{|k|^{'}}_{\h_{\lambda}^{\bot},\mathfrak{q}}(F_p^{(i)})+\sum_{\substack{i,k\\i+p=\gamma-1\\i-k=b-1}}d^{|k|^{''}}_{\h_{\lambda}^{\bot},\mathfrak{q}}(F_p^{(i)})\right)\right)=\]
\[=\sum_{\gamma}t^{-\gamma+1}\left(\sum_b\left(\sum_{\substack{i,k\\i+p=\gamma\\i-k=b}} d^{|k|^{'}}_{\h_{\lambda}^{\bot},\mathfrak{q}}(F_p^{(i)})+\sum_{\substack{i,k\\i+p=\gamma\\i-k=b}}d^{|k|^{''}}_{\h_{\lambda}^{\bot},\mathfrak{q}}(F_p^{(i-1)})\right)\right).\]
The summation index in equations (\ref{firsteq2}) gives $i+p=a+b$. So for $a+b=\gamma$ we get

\begin{equation}\label{second2}
\sum_{\substack{i,k\\i+p=\gamma\\i-k=b}} d^{|k|^{'}}_{\h_{\lambda}^{\bot},\mathfrak{q}}(F_p^{(i)})+\sum_{\substack{i,k\\i+p=\gamma\\i-k=b}}d^{|k|^{''}}_{\h_{\lambda}^{\bot},\mathfrak{q}}(F_p^{(i-1)})=0,
\end{equation}

and so $G_t\in H^0(\h_{t\lambda}^{\bot},d_{\h_{t\lambda}^{\bot},\mathfrak{q}})$, $\forall t\in\mathbb{R}^{\ast}$.  $\diamond$

Thus if $F\in H^0(\h_{\lambda}^{\bot},d_{\h_{\lambda}^{\bot},\mathfrak{q}})$, then $F\in \overline{J}\left(H^0_{(\epsilon=1)}(\h_{\lambda}^{\bot},d^{(\epsilon=1)}_{\h_{\lambda}^{\bot},\mathfrak{q}})\right)$ iff there is a family $(F_p)_{1\geq p\geq n}\in S(\mathfrak{q})$ such that $F_{(t)}:=\sum_pt^pF_p$ satisfies the conditions  $\forall t\in\mathbb{R}^{\ast},\;\;F_{(t)}\in H^0(\h_{t\lambda}^{\bot},d_{\h_{t\lambda}^{\bot},\mathfrak{q}})$ and $F_{(t=1)}=F$.

\textbf{6.3 Deformation algebra.} Let $\g$ be a Lie algebra. Consider a supplementary variable $T$ such that $[T,\g]=0$ and set $\g_T=\g\oplus <T>$ and $\h_T=\h\oplus <T>$ such that $\dim(\g_T)=\dim(\g)+1$. Set also $U(\g_T)$ to be the U.E.A of $\g_T$ and $U(\g_T)\h^T_{\lambda}$ to be the ideal of $U(\g_T)$ generated by  $\h_{\lambda}^T=<H+T\lambda(H),\;H\in\h>$. Set $\mathbf{D}_T(\g,\h,\lambda):=(U(\g_T)/U(\g_T)\h^T_{\lambda})^{\h_T}$. We use the notation $\overline{\beta}_t:\;S(\mathfrak{q})\stackrel{\sim}{\longrightarrow} U(\g)/U(\g)\h_{t\lambda}$ to denote this vector space isomorphism through symmetrization.
\newtheorem{dkr}[con]{Definition}
\begin{dkr}
Let $\g$ be a Lie algebra and $\h\subset\g$ a subalgebra and $\lambda$ a character of $\h$. Fix $\mathfrak{q}$ a supplementary of $\h$. We will say that $t\longrightarrow u_t\in U(\g)/U(\g)\h_{t\lambda}$ is a polynomial family in $t$ iff $t\longrightarrow \overline{\beta}_t^{-1}(u_t)\in S(\mathfrak{q})$ is a polynomial family in $t$. We also denote as $\mathcal{P}_{(t)}\left((U(\g)/U(\g)\h_{t\lambda})^{\h}\right)$ the polynomial in $t$ families  $t\longrightarrow u_t\in\left(U(\g)/U(\g)\h_{t\lambda}\right)^{\h}.$  Equivalently for such a family we will write $u_t\in\mathcal{P}_{(t)}\left((U(\g)/U(\g)\h_{t\lambda})^{\h}\right)$.
\end{dkr}

\newtheorem{stre}[con]{Lemma}
\begin{stre}
The object $\mathcal{P}_{(t)}\left((U(\g)/U(\g)\h_{t\lambda})^{\h}\right)$ is an algebra.
\end{stre}
\textit{Proof.} For polynomial in $t$ families $u_t$ and $v_t$, $u_t+v_t$ is in $(U(\g)/U(\g)\h_{t\lambda})^{\h}$ and polynomial in $t$. Let us verify that if $\overline{\beta}_t^{-1}(u_t), \overline{\beta}_t^{-1}(v_t)$ are polynomial families in $t$, then $t\longrightarrow \overline{\beta}_t^{-1}(u_tv_t)$ is polynomial in $t$ too. Let $\{X_1,\ldots,X_p\}$ be a basis of $\mathfrak{q}$ and $\{H_1,\ldots,H_r\}$ be a basis of $H$ such that $\{H_2,\ldots,H_r\}$ is a basis of $ker \lambda\cap \h$. For multiindices $\alpha,\nu,\gamma,\in \mathbb{N}^p$, $\delta\in \mathbb{N}^r$ we use the notation $X^{\alpha}=X_1^{\alpha_1}\cdots X_p^{\alpha_p},\;X^{\gamma}=X_1^{\gamma_1}\cdots X_p^{\gamma_p},\;H^{\delta}=H_1^{\delta_1}\cdots H_r^{\delta_r}$. With the assumptions on the basis of $\h$, we have that $H^{\delta}\in U(\g)\h_{t\lambda}$ if $(\delta_2,\ldots,\delta_r)\neq (0,\ldots,0)$. Then 
\[\beta(X^{\alpha})\beta(X^{\gamma})=\sum_{\alpha,\gamma,\nu}c_{\alpha\gamma}^{\nu}\beta(X^{\nu})+\sum_{\alpha,\gamma,\nu,\delta}c_{\alpha\gamma}^{\nu\delta}\beta(X^{\nu})H^{\delta}\equiv \sum_{\alpha,\gamma,\nu}c_{\alpha\gamma}^{\nu}\beta(X^{\nu})\cdot (-t\lambda(H_1))^{\delta_1}\;mod[U(\g)\h_{t\lambda}].\]
Using the decomposition $U(\g)=\beta(S(\mathfrak{q}))\oplus U(\g)\h_{\lambda}$ of $U(\g)$, we get that the element $R=\overline{\beta}_t^{-1}\left(\beta(X^{\alpha})\beta(X^{\gamma})\right)=\sum_{\alpha,\gamma,\nu}c_{\alpha\gamma}^{\nu}X^{\nu}+\sum_{\alpha,\gamma,\nu}c_{\alpha\gamma}^{\nu\delta_1}X^{\nu}(-t\lambda(H_1))^{\delta_1}$ is polynomial in $t$. Now write $\overline{\beta}^{-1}_t(u_t)=\sum_{\alpha}p_{\alpha}(t)X^{\alpha}$ and $\overline{\beta}^{-1}_t(v_t)=\sum_{\gamma}q_{\gamma}(t)X^{\gamma}$ as polynomial functions. Then $\overline{\beta}^{-1}_t(u_tv_t)=\sum_{\alpha,\gamma}p_{\alpha}(t)q_{\gamma}(t)\overline{\beta}^{-1}_t\left(\beta_t(X^{\alpha})\beta_t(X^{\gamma})\right)$ is a polynomial family in $t$. $\diamond$

Let $<T-t>$ denote the corresponding ideal of $U(\g_T)$ and $\pi_t:\;U(\g_T)\longrightarrow U(\g_T)/<T-t>$ denote the canonical projection. Let $S_T(\g)$ be the symmetric algebra of $\g_T$. There is a surjective algebra map $e_t:\;U(\g_T)\longrightarrow U(\g)$ defined by $T\mapsto t$ and $\forall X\in\g,\; X\mapsto X$. The kernel of $e_t$ is $U(\g_T)<T-t>$ and we get an isomorphism of algebras $\left(U(\g_T)/<T-t>\right)\simeq U(\g)$. We write $\pi_t\left( (U(\g_T)/U(\g_T)\h^T_{\lambda})^{\h_T}\right)=\left( (U(\g_T)/U(\g_T)\h^T_{\lambda})^{\h_T}/<T-t>\right)$. So by evaluation at $T=t$, we get an injective map of algebras 
\[e_{(T=t)}:\;\left( (U(\g_T)/U(\g_T)\h^T_{\lambda})^{\h_T}/<T-t>\right)\hookrightarrow (U(\g)/U(\g)\h_{t\lambda})^{\h}.\]

\newtheorem{okra2}[con]{Theorem}
\begin{okra2}\label{ha}
Let $\g$ be a Lie algebra, and $\h,\lambda,\mathfrak{q}$ as usual. Let $t\longrightarrow u_t\in (U(\g)/U(\g)\h_{t\lambda})^{\h}$ be a polynomial family in $t$. Then there is $u_T\in (U(\g_T)/U(\g_T)\h^T_{\lambda})^{\h_T}$ such that $e_t(u_T)=u_t$.
\end{okra2}
\textit{Proof.} Choose a PBW basis $\{Q_1,\ldots,Q_p,H_1,\ldots H_m,X\}$ for $\g=\mathfrak{q}\oplus\h_0\oplus<X>$ where $X\in\h$ and $\h_0$ be such that $\g=\mathfrak{q}\oplus\h_0\oplus<X>$ and $\lambda(X)=1$. The PBW theorem gives a supplementary of $U(\g)\h_{t\lambda}$ in $ U(\g)$. In fact considering the symmetrization map $\overline{\beta}_t:\;S(\mathfrak{q})\longrightarrow U(\g)/U(\g)\h_{t\lambda}$ and $Q^{\alpha}=Q_1^{\alpha_1}\cdots Q_p^{\alpha_p}$ a basis of $S(\mathfrak{q})$, then for $F(t)\in \mathcal{P}_{(t)}\left((U(\g)/U(\g)\h_{t\lambda})^{\h}\right)$, $\exists\; (p_\alpha(t))$ a family of polynomials such that $F(t)$ can be written as $F(t)=\sum_\alpha p_{\alpha}(t)\overline{\beta}_t(Q^{\alpha})$. The condition for $F(t)$ to be $\h-$ invariant is 
\begin{equation}\label{invar}
[H_i,F(t)]\in U(\g)\h_{t\lambda},\;\;\forall i.
\end{equation}
For $Q_1,\ldots,Q_p$ a PBW basis, $\sigma,\alpha,\gamma,\delta$ multiindices, $Q^\sigma:=Q_1^{\sigma_1}\cdots Q_p^{\sigma_p}$ and $\beta:\;S(\g)\longrightarrow U(\g)$ the ordinary symmetrization, we write that $\forall i,\;\;[H_i,\beta(Q^{\alpha})]=\sum c_{\sigma\gamma\delta}^{(i)(\alpha)}Q^\sigma H^\gamma X^\delta.$ Since $F(t)=\sum_\alpha p_{\alpha}(t)\overline{\beta}_t(Q^{\alpha})$, condition (\ref{invar}) is equivalent to 
\begin{equation}\label{act}
\sum_\alpha p_{\alpha}(t)\sum_{i,\sigma,\delta} c_{\sigma 0\delta}^{(i)(\alpha)}Q^{\sigma}(-t)^{\delta}=0.
\end{equation}
Let now $\beta_T:\;S_T(\g)\longrightarrow U(\g_T)$ be the corresponding symmetrization map and consider the element $u_T:=\sum_{\alpha}P_{\alpha}(T)\beta_T(Q^{\alpha})$. The condition for $u_T$ to be $\h_T-$ invariant is
\begin{equation}\label{sect}
\sum_{\alpha}P_{\alpha}(T)\sum_{i,\sigma,\delta} c_{\sigma 0\delta}^{(i)(\alpha)}Q^{\sigma}(-T)^{\delta}=0.
\end{equation}
Clearly the equations (\ref{act}) and (\ref{sect}) are equivalent and $u_T$ constructed above satisfies $\pi_t(u_T)=u_t$.

\newtheorem{corm}[con]{Corollary}
\begin{corm}\label{hi}
In the specialized version $t=1$ of theorem \ref{ha} we have for the algebra $\mathbf{D}_{(T=1)}(\g,\h,\lambda):=  (U(\g_T)/U(\g_T)\h^T_{\lambda})^{\h_T}/<T-~1>$ that 
\[\mathbf{D}_{(T=1)}(\g,\h,\lambda)\hookrightarrow (U(\g)/U(\g)\h_{\lambda})^{\h},\]
and in particular, the elements of $\mathbf{D}_{(T=1)}(\g,\h,\lambda)$ correspond to elements who are the value at $t=1$ of elements in $\mathcal{P}_{(t)}\left((U(\g)/U(\g)\h_{t\lambda})^{\h}\right)$: If $u\in \mathbf{D}_{(T=1)}(\g,\h,\lambda)$ then there is an element $u_T\in U(\g_T)$ such that $u=\pi_{(T=1)}(u_T)$. Moreover, the element $u_t:=e_{(T=t)}(u_T)\in (U(\g)/U(\g)\h_{t\lambda})^{\h}$ defines a polynomial family in $t$, that is $u_t\in \mathcal{P}_{(t)}\left((U(\g)/U(\g)\h_{t\lambda})^{\h}\right)$.
\end{corm}
We denote abusively and for presentation reasons as $\mathcal{P}_{(t=1)}\left((U(\g)/U(\g)\h_{\lambda})^{\h}\right)$ the values at $t=1$ of elements $u_t\in \mathcal{P}_{(t)}\left((U(\g)/U(\g)\h_{t\lambda})^{\h}\right)$. Recall from Lemma \ref{jbar} that  $H^0_{(\epsilon=1)}(\h_{\lambda}^{\bot},d^{(\epsilon=1)}_{\h_{\lambda}^{\bot},\mathfrak{q}})\hookrightarrow H^0(\h_{\lambda}^{\bot},d_{\h_{\lambda}^{\bot},\mathfrak{q}})$ and that from $\mathcal{x}$ 6.1, $H^0(\h_{\lambda}^{\bot},d_{\h_{\lambda}^{\bot},\mathfrak{q}})\hookrightarrow (U(\g)/U(\g)\h_{\lambda+\rho})^{\h}$. Thus $H^0_{(\epsilon=1)}(\h_{\lambda}^{\bot},d^{(\epsilon=1)}_{\h_{\lambda}^{\bot},\mathfrak{q}})\hookrightarrow(U(\g)/U(\g)\h_{\lambda})^{\h}$. We denote as $\mathfrak{i}_{(\epsilon=1)}$ the injective map $\mathfrak{i}_{(\epsilon=1)}:\;H^0_{(\epsilon=1)}(\h_{\lambda}^{\bot},d^{(\epsilon=1)}_{\h_{\lambda}^{\bot},\mathfrak{q}})\hookrightarrow(U(\g)/U(\g)\h_{\lambda+\rho})^{\h}$.

\newtheorem{okri1}[con]{Proposition}
\begin{okri1}\label{ho}
Let $\h_{\lambda}^T=<H+\lambda(H)T,\;H\in\h>\subset \g_T$ and $\mathfrak{q}$ a supplementary of $\h$ in $\g$. Let also $t\mapsto F_t \in H^0(\h_{t\lambda}^{\bot},d_{\h_{t\lambda}^{\bot},\mathfrak{q}})$ be a polynomial family in $t$. Then there is $F_T\in  H^0((\h_{\lambda}^T)^{\bot},d_{(\h_{\lambda}^T)^{\bot},\mathfrak{q}})$ such that $e_t(F_T)=F_t$.
\end{okri1}
\textit{Proof.}  Let $F_t=\sum_kt^kP_k$ with $\forall k,\;P_k\in S(\mathfrak{q})$ and $F_t\in  H^0(\h_{t\lambda}^{\bot},d_{\h_{t\lambda}^{\bot},\mathfrak{q}})$, that is $d_{\h_{t\lambda}^{\bot},\mathfrak{q}}(F_t)=0$. Analysing the differential $d_{\h_{t\lambda}^{\bot},\mathfrak{q}}$ in terms of $\mathrm{deg}_{\mathfrak{q}}$ and $F_t$ in homogeneous polynomials $P_k^{(i)}$ as in the proof of Theorem \ref{difo}, we conclude that similarly to (\ref{firsteq1}), the equation $d_{\h_{t\lambda}^{\bot},\mathfrak{q}}(F_t)=0$ is equivalent to
\begin{equation}\label{third eq}
\sum_{\substack{i,s\\i-s=b}}d^{|s|^{'}}_{\h_{\lambda}^{\bot},\mathfrak{q}}(P_k^{(i)})+d^{|s|^{''}}_{\h_{\lambda}^{\bot},\mathfrak{q}}(P_{k+1}^{(i-1)})=0.
\end{equation}
Now since $T$ is a central variable, we have $\forall H\in\h,\;[H+\lambda(H)T,\cdot]=[H,\cdot]$, and so the differential $d_{(\h_{\lambda}^T)^{\bot},\mathfrak{q}}$ contains exactly the same graphs as $d_{\h^{\bot},\mathfrak{q}}$. 
Decomposing $d_{(\h_{\lambda}^T)^{\bot},\mathfrak{q}}$ as $d_{(\h_{\lambda}^T)^{\bot},\mathfrak{q}}=d_{\h^{\bot},\mathfrak{q}}=\sum_k d^{|s|}_{\h^{\bot},\mathfrak{q}}$ with $\mathrm{deg}_{\mathfrak{q}}(d^{|s|}_{\h^{\bot},\mathfrak{q}})=-s+1$, the equation (\ref{third eq}) implies that $F_T:=\sum_k T^{\ast^k}P_k$ satisfies the equation $d_{(\h_{\lambda}^T)^{\bot},\mathfrak{q}}(F_T)=0$ and thus $F_T\in H^0((\h_{\lambda}^T)^{\bot},d_{(\h_{\lambda}^T)^{\bot},\mathfrak{q}})$ and $e_t(F_T)=F_t$.
\newtheorem{coro}[con]{Corollary}
\begin{coro}\label{he}
Proposition \ref{ho} specialized at $T=1$ says that $H^0((\h_{\lambda}^T)^{\bot},d_{(\h_{\lambda}^T)^{\bot},\mathfrak{q}})/<T-1>$ is the values at $t=1$ of polynomial in $t$ families $t\longrightarrow F_t\in H^0(\h_{t\lambda}^{\bot},d_{\h_{t\lambda}^{\bot},\mathfrak{q}})$.
\end{coro}

\newtheorem{okri}[con]{Theorem}
\begin{okri}\label{final}
The specialized algebras $\mathbf{D}_{(T=1)}(\g,\h,\lambda+\rho):= \left((U(\g_T)/U(\g_T)\h^T_{\lambda+\rho})^{\h_T}/<T-1>\right)$ and \newline $H^0_{(\epsilon=1)}(\h_{\lambda}^{\bot},d^{(\epsilon=1)}_{\h_{\lambda}^{\bot},\mathfrak{q}}):=\left(H^0_{(\epsilon)}(\h_{\lambda}^{\bot},d^{(\epsilon)}_{\h_{\lambda}^{\bot},\mathfrak{q}})/<\epsilon-1>\right)$ are isomorphic.
\end{okri}
\textit{Proof.} In the vector space case that we consider it is
\begin{equation}\label{isoT}
(U(\g_T)/U(\g_T)\h^T_{\lambda+\rho})^{\h_T}\simeq H^0((\h_{\lambda}^T)^{\bot},d_{(\h_{\lambda}^T)^{\bot},\mathfrak{q}}).
\end{equation}
Now we specialize (\ref{isoT}) at $T=1$ and see what we get from both sides. At the left hand side we have $\mathbf{D}_{(T=1)}(\g,\h,\lambda+\rho)=\left((U(\g_T)/U(\g_T)\h^T_{\lambda+\rho})^{\h_T}/<T-1>\right)$ which by Theorem \ref{ha} and Corollary \ref{hi} are the values at $t=1$ of polynomial in $t$ families $t\longrightarrow u_t\in (U(\g)/U(\g)\h_{t(\lambda+\rho)})^{\h}$. At the right hand side of (\ref{isoT}) specialized at $T=1$, we have $H^0((\h_{\lambda}^T)^{\bot},d_{(\h_{\lambda}^T)^{\bot},\mathfrak{q}})/<T-1>$ which by Corollary \ref{he} is the values at $t=1$ of the polynomial in $t$ families $t\longrightarrow F_t\in H^0(\h_{t\lambda}^{\bot},d_{\h_{t\lambda}^{\bot},\mathfrak{q}})$.  By Theorem \ref{difo},  $H^0_{(\epsilon=1)}(\h_{\lambda}^{\bot},d^{(\epsilon=1)}_{\h_{\lambda}^{\bot},\mathfrak{q}})$ is also the values at $t=1$ of polynomial in $t$ families $t\longrightarrow F_t\in H^0(\h_{t\lambda}^{\bot},d_{\h_{t\lambda}^{\bot},\mathfrak{q}})$. We have thus shown that
\[H^0_{(\epsilon=1)}(\h_{\lambda}^{\bot},d^{(\epsilon=1)}_{\h_{\lambda}^{\bot},\mathfrak{q}})\simeq\mathbf{D}_{(T=1)}(\g,\h,\lambda+\rho)\hookrightarrow (U(\g)/U(\g)\h_{\lambda+\rho})^{\h}. \;\;\diamond\]

\end{document}